\documentclass{amsart}
\usepackage{amssymb, amsmath, amscd, latexsym, mathrsfs, eufrak, lscape, tabls}
\input xy
\xyoption{all}
\newtheorem{thm}{Theorem}[section]
\newtheorem{cor}[thm]{Corollary}
\newtheorem{lem}[thm]{Lemma}

\newtheorem{prop}[thm]{Proposition}
\newtheorem{defn}[thm]{Definition}

\newtheorem{rem}[thm]{Remark}

\newtheorem{claim}[thm]{Claim}
\newtheorem{ex}{Example}[section]
\numberwithin{equation}{section}

\newtheorem{bst}{Bisymmetric Triple}[section]
\newtheorem{symp}{Symmetric Pair}[section]

\newenvironment{pr}{\textit{Proof: }}

\newcommand{\sg}{\sigma}

\newcommand{\al}{\alpha}
\newcommand{\be}{\beta}
\newcommand{\la}{\lambda}

\newcommand{\ga}{\gamma}

\newcommand{\de}{\delta}

\newcommand{\iso}{\cong}

\newcommand{\li}{\medskip}

\newcommand{\reals}{\mathbb{R}}

\newcommand{\complex}{\mathbb{C}}

\newcommand{\R}{\mathcal{R}}

\newcommand{\mfe}{\mathfrak{e}}
\newcommand{\mff}{\mathfrak{f}}
\newcommand{\mfg}{\mathfrak{g}}
\newcommand{\mfh}{\mathfrak{h}}
\newcommand{\mfk}{\mathfrak{k}}
\newcommand{\mfl}{\mathfrak{l}}
\newcommand{\mfm}{\mathfrak{m}}
\newcommand{\mfn}{\mathfrak{n}}

\newcommand{\mfp}{\mathfrak{p}}

\newcommand{\mfu}{\mathfrak{u}}

\newcommand{\mfsu}{\mathfrak{s}\mathfrak{u}}
\newcommand{\mfso}{\mathfrak{s}\mathfrak{o}}
\newcommand{\mfsp}{\mathfrak{s}\mathfrak{p}}

\def\bar{\begin{array}}
\def\ear{\end{array}}
\def\sbar{\begin{subarray}}
\def\sear{\end{subarray}}
\def\beq{\begin{equation}}
\def\eeq{\end{equation} }
\def\beqar{\begin{eqnarray}}
\def\eeqar{\end{eqnarray}}
\def\bal{\begin{align}}
\def\eal{\end{align}}
\def\bfig{\begin{figure}}
\def\efig{\end{figure}}
\def\bc{\begin{center}}
\def\ec{\end{center}}
\def\btab{\begin{table}}
\def\etab{\end{table}}
\def\bland{\begin{landscape}}
\def\eland{\end{landscape}}

\def\bproof{\begin{pr}}
\def\eproof{\end{pr}}
\def\bprop{\begin{prop}}
\def\eprop{\end{prop}}
\def\bthm{\begin{thm}}
\def\ethm{\end{thm}}
\def\blem{\begin{lem}}
\def\elem{\end{lem}}
\def\brem{\begin{rem}}
\def\erem{\end{rem}}
\def\bcor{\begin{cor}}
\def\ecor{\end{cor}}
\def\bex{\begin{ex}}
\def\eex{\end{ex}}
\def\bdfn{\begin{defn}}
\def\edfn{\end{defn}}
\def\bclaim{\begin{claim}}
\def\eclaim{\end{claim}}
\def\bbst{\begin{bst}}
\def\ebst{\end{bst}}
\def\bsymp{\begin{symp}}
\def\esymp{\end{symp}}

\newlength{\myVSpace}
\newcommand\xstrut{\raisebox{-.9\myVSpace}{\rule{0pt}{\myVSpace}}}

\makeatletter
\def\table{\@ifnextchar[{\table@i}{\table@i[\fps@table]}}
\def\table@i[#1]{\@float{table}[#1]\footnotesize}
\makeatother

\begin{document}

\title[]
{Einstein Homogeneous Bisymmetric Fibrations}

\author{F\'{a}tima Ara\'{u}jo}

\address{School of Mathematics, The University of Edinburgh, JCMB, The Kings Buildings, Edinburgh, EH9 3JZ}

\email{m.d.f.araujo@sms.ed.ac.uk}


\date{}%
\maketitle

\begin{abstract}
We consider a homogeneous fibration $G/L \rightarrow G/K$, with symmetric fiber and base, where $G$ is a compact connected semisimple Lie group and $L$ has maximal rank in $G$. We suppose the base space $G/K$ is isotropy irreducible and the fiber $K/L$ is simply connected. We investigate the existence of $G$-invariant Einstein metrics on $G/L$ such that the natural projection onto $G/K$ is a Riemannian submersion with totally geodesic fibers. These spaces are divided in two types: the fiber $K/L$ is isotropy irreducible or is the product of two irreducible symmetric spaces. We classify all the $G$-invariant Einstein metrics with totally geodesic fibers for the first type. For the second type, we classify all these metrics when $G$ is an exceptional Lie group. If $G$ is a classical Lie group we classify all such metrics which are the orthogonal sum of the normal metrics on the fiber and on the base or such that the restriction to the fiber is also Einstein.
\end{abstract}

\maketitle

\section{Introduction}\label{intro}

\addtolength{\myVSpace}{0.1cm}

A Riemannian manifold $(M,g)$ is said to be Einstein if its Ricci curvature satisfies an equation of the form $Ric=Eg$, for some constant $E$. This equation is a system of second order partial differential equations, which is in general unmanageable. Fully general results are not known, but many results of existence and classification of Einstein metrics are known for some special manifolds. Example of this are the K\"{a}hler-Einstein (\cite{Yau}) and the Sasakian-Einstein manifolds (\cite{BG}). For a homogeneous space the Einstein equation is a system of algebraic equations, which is an easier problem than its general version. Due to this, most of the known examples of existence or non-existence of Einstein metrics are homogeneous spaces. For example, every isotropy irreducible space is clearly an Einstein manifold and recently Einstein metrics on homogeneous spaces with exactly two isotropy summands were classified by Dickinson and Kerr (\cite{DK}). Einstein metrics on spheres and projective spaces were classified by Ziller (\cite{Zi2}) and Einstein normal homogeneous manifolds were classified by Wang and Ziller (\cite{WZ}). It is known that every compact simply connected homogeneous manifold with dimension less or equal to $11$ admits a homogeneous Einstein metric (\cite{ADF}, \cite{Be}, \cite{BK}, \cite{CRW}, \cite{Je1}, \cite{NR2}). There are examples in dimension $12$ which do not admit an Einstein homogeneous metric (\cite{BK}, \cite{WZ2}). For a survey on results for Einstein manifolds, see \cite{Wa}. Riemannian submersions have also been used to construct new Eintein metrics. We recall the work of Jensen on principal  fibers bundles (\cite{Je2}) and the work of Wang and Ziller on principal torus bundles (\cite{WZ3}). We use some results obtained by the author in \cite{Fa2} about Einstein homogeneous fibrations to investigate the existence of Einstein metrics on the total space of fibrations whose fiber and base are symmetric spaces.

\li

Let $G$ be a compact connected semisimple Lie group and $L\varsubsetneq K\varsubsetneq G$
connected closed non-trivial subgroups, such that $G/K$ is an irreducible symmetric space, $K/L$ is a simply-connected symmetric space, not necessarily irreducible, and $L$ has maximal rank. We consider the homogeneous fibration

\beq\label{bsfibration} M = G/L \rightarrow G/K = N \textrm{ with fiber } F = K/L.\eeq

Throughout, we call a fibration as above a \textbf{bisymmetric fibration}. We consider on $M$ the class of metrics $g_M$ such that the natural projection $M \ni aL\mapsto aK\in N,\,a\in G,$ is a Riemannian submersion with totally geodesic fibers. Throughout a metric with this property is called an \textbf{adapted} metric (see \cite{Fa2}). The aim of this paper is to classify $G$-invariant Einstein adapted metrics on bisymmetric fibrations.

\li

Let $\mfg$, $\mfk$ and $\mfl$ denote the Lie algebras of $G$, $K$ and $L$, respectively. Let $\Phi$ be the Killing form of $G$ and consider the $Ad\,G$-invariant symmetric bilinear form given by $B=-\Phi$. Since $G$ is compact and semisimple, $B$ is positive definite. We consider a $B$-orthogonal decomposition of $\mfg$ given by

\beq\mfg=\mfl\oplus\mfm = \mfl\oplus\mfp\oplus\mfn,\eeq

where $\mfg=\mfl\oplus\mfm$, $\mfg=\mfk\oplus\mfn$ and $\mfk=\mfl\oplus\mfp$ are reductive
decompositions for $M$, $N$ and $F$, respectively. We consider a $B$-orthogonal decomposition $\mfp=\mfp_1\oplus\ldots\oplus\mfp_s$ of $\mfp$ into irreducible $Ad\,L$-modules. We note that $\mfn$ is an irreducible $Ad\,K$-module, but not necessarily $Ad\,L$-irreducible.

We denote by $B_{\mathfrak{q}}$ the restriction of $B$ to some subspace $\mathfrak{q}\subset \mfg $. A metric on $M$ defined by an $Ad\,L$-invariant Euclidean product on $\mfm$ of the form

\beq\label{mdefini}g_{\mfm}=\left(\oplus_{a=1}^s\la_aB_{\mfp_a}\right)\oplus\mu B_{\mfn},\,\la_a,\,\mu>0\eeq

is an adapted metric. Throughout, we assume the following hypothesis:

\beq\label{submodulesHyp}\bar{l}\mfp_1,\ldots,\mfp_s \textrm{ are pairwise
inequivalent irreducible }Ad\,L\textrm{-submodules;}\\
\mfp \textrm{ and } \mfn \textrm{ do not contain equivalent }
Ad\,L\textrm{-submodules.} \ear\eeq

Under the hypothesis (\ref{submodulesHyp}), any adapted metric on $M$ is determined by an $Ad\,L$-invariant Euclidean product given by (\ref{mdefini}) (see \cite{Fa2}). We denote an adapted metric on $M$ by $g_M$, its restriction to the fiber by $g_F$ and its projection onto the base space by $g_N$.

\li

For a bisymmetric fibration as in (\ref{bsfibration}), we call $(\mfg,\mfk,\mfl)$ a \textbf{bisymmetric triple} of maximal rank. The triple $(\mfg,\mfk,\mfl)$ is said to be irreducible is $\mfg$ is simple. Clearly, there is a one-to-one correspondence between bisymmetric fibrations, up to cover, and bisymmetric triples. All the bisymmetric triples considered in this paper are irreducible and such that $\mfl$ has maximal rank.

\li

The classification of isotropy irreducible symmetric spaces is very well known and can be found in \cite{He}. By using this classification we obtain a list of all possible triples $(\mfg,\mfk,\mfl)$ such that $\mfl$ and
$\mfk$ are subalgebras of maximal rank of $\mfg$ and $(\mfg,\mfk)$
and $(\mfk,\mfl)$ are symmetric pairs. In particular, any irreducible bisymmetric triple of maximal rank is of Type I or II, where $(\mfg,\mfk,\mfl)$ is said to be of \textbf{Type I} if $F=K/L$ is an isotropy irreducible symmetric space and of \textbf{Type II} if $F=K/L$ is the direct product of two isotropy irreducible symmetric spaces. A bisymmetric fibration $M=G/L\rightarrow G/K= N$ is said to be of Type I or II if the corresponding bisymmetric triple
$(\mfg,\mfk,\mfl)$ is either of Type I or II, respectively. We present a list of all irreducible bisymmetric triples of maximal rank. More precisely we will show the following result in Section \ref{bsSection}:

\li

\blem \label{bsclass} An irreducible bisymmetric triple $(\mfg,\mfk,\mfl)$ of maximal rank is either of Type I
or II. Moreover, all such bisymmetric triples are those in Tables and \ref{eigIexc}, \ref{eigIclass},
\ref{eigIIexc} and \ref{eigIIclass}.\elem

\li

We classify all the Einstein adapted metrics for bisymmetric fibrations of type I. Since the fiber $F$ for a bisymmetric fibration of type I is an isotropy irreducible symmetric space, any adapted metric $g_M$ is determined by an $Ad\,L$-invariant Euclidean product of the form

\beq\label{mdefbiini}g_{\mfm}=\la B_{\mfp}\oplus\mu B_{\mfn}.\eeq

A metric of this form is called a \textbf{binormal} metric (see \cite{Fa2}). We remark that the fact that the base space $N$ and the base $F$ are isotropy irreducible does not imply that $M$ has only two isotropy submodules. Indeed, the horizontal subspace $\mfn$ is $Ad\,K$-irreducible but is not in general $Ad\,L$-irreducible. The particular case of existence of $G$-invariant Einstein metrics when $M$ has exactly two irreducible isotropy subspaces was studied by McKenzie Y. Wang and Ziller, under some assumptions, in \cite{WZ2} and, more recently and in full generality, by W. Dickinson and M. M. Kerr in \cite{DK}, who classified all such metrics. For more details on the decomposition of the isotropy representation for irreducible bisymmetric fibrations of maximal rank see \cite{Fa1}. In Section \ref{bsSection} we prove the following result.

\bthm \label{mtypeI} The bisymmetric fibrations $M=G/L\rightarrow G/K$ of Type I such that $M$ admits an Einstein adapted metric are those in Tables
\ref{mIexc} and \ref{mIclass}. For each case there are exactly two
Einstein adapted metrics. Furthermore, these Einstein metrics
are, up to homothety, given by

$$g_{\mfm}=B_{\mfp}\oplus XB_{\mfn},$$ where the two distinct values of $X$ are indicated in the tables mentioned above. In all the cases, $g_N$ and
$g_F$ are also Einstein. \ethm

For bisymmetric fibrations of type II the vertical space $\mfp$ decomposes into two irreducible $Ad\,L$-submodules. Hence, an adapted metric is not necessarily binormal (see (\ref{mdefbiini})). Moreover, there may exist adapted metrics whose restriction to the fiber is not Einstein. Theorems \ref{biII} and \ref{gfeinII} classify all the bisymmetric fibrations of type II which admit an Einstein binormal metric or an Einstein adapted metric whose restriction to the fiber is also Einstein. These results are proved in Section \ref{bsSection}.

\bthm \label{biII} The bisymmetric fibrations $M=G/L\rightarrow G/K$ of Type II such that $M$ admits an Einstein binormal metric are those listed
in Table \ref{bimII}. Furthermore, the binormal Einstein metrics
are, up to homothety, given by

$$g_{\mfm}=B_{\mfp}\oplus XB_{\mfn},$$ where
$X$ is indicated in Table \ref{bimII}. In all the cases, $g_N$ and
$g_F$ are also Einstein.

In particular, if $M$ admits an Einstein binormal metric, then $G$ is a classical Lie group.\ethm

\bthm \label{gfeinII}The bisymmetric fibrations $M=G/L\rightarrow G/K$ of Type II such that $M$ admits an
Einstein adapted metric $g_M$ whose restriction $g_F$ to the fiber $F=K/L$ is also Einstein are

(i) those with an Einstein binormal metric, given by Theorem \ref{biII} and Table
\ref{bimII};

(ii) the fibration given by

$$(\mfsu_{2(l+s)},\mfsu_{2l}\oplus\mfsu_{2s}\oplus\reals,\mfsu_l\oplus\mfsu_l\oplus\mfsu_s\oplus\mfsu_s\oplus\reals^3),$$

whose Einstein adapted metric is, up to homothety, given by

$$g_{\mfm}=\frac{2l}{l+s}B_{\mfp_1}\oplus\frac{2s}{l+s}B_{\mfp_2}\oplus B_{\mfn}.$$

The metric in (ii) is binormal if and only if $l=s$.

Furthermore, if $M$ admits an Einstein adapted metric whose restriction to the fiber is also Einstein, then $G$ is a classical Lie group.\ethm

\li

According to Theorems \ref{biII} and \ref{gfeinII}, if $G$ is an exceptional Lie group it does not admit an Einstein binormal metric or an Einstein adapted metric which restricts to an Einstein metric on the fiber. However, in the exceptional case it is possible to classify all the Einstein adapted metrics. In Section \ref{bsSection} we show the following:

\bthm  \label{genIIexc} The only bisymmetric fibrations $M=G/L\rightarrow G/K$ of Type II, for an exceptional Lie group $G$, which admit an Einstein adapted metric are those listed in Table \ref{tabgenII}. The Einstein adapted metrics are, up to homothety, given by,

$$g_{\mfm}=\frac{1}{X_1}B_{\mfp_1}\oplus\frac{1}{X_2}B_{\mfp_2}\oplus B_{\mfn}$$

and approximations for $X_1$, $X_2$ are given in Table \ref{tabgenII}. These metrics are not binormal and the restriction to the fiber is not Einstein.\ethm

Unlike the exceptional case, if $G$ is a classical Lie group, due to the complexity of the Einstein equations, we do not present a full classification of Einstein adapted metrics for type II. Still it is possible to classify all the Einstein adapted metrics for type II if the two eigenvalues $\ga_1$, $\ga_2$ of the Casimir operator of $\mfk$ on $\mfp_1$, $\mfp_2$, respectively, satisfy the condition $\ga_2=\ga_1$ or $\ga_2=1-\ga_1$. This is the case when the restriction $g_F$ of an Einstein adapted metric $g_M$ to the fiber is Einstein (see Section \ref{sectionsymF}). In particular, this implies that if $M$ admits an Einstein adapted metric $g_M$ such that $g_F$ is also Einstein, then we can classify all the other Einstein adapted metrics on $M$.

We recall that if $U$ is a vector subspace of $\mfg$, the Casimir operator of $U$ is the operator

\beq \label{casdef}C_U=\sum_i(ad_{u_i})^2\in \mfg\mfl(\mfg),\eeq

where $\{u_i\}_i$ is an orthonormal basis of $U$ with respect to $\Phi$. The result is as follows:

\bthm \label{genIIall} Let $M=G/L\rightarrow G/K$ be a bisymmetric fibration
of Type II. Let $\ga_1$, $\ga_2$ be the eigenvalues of the Casimir operator $C_{\mfk}$ on $\mfp_1$, $\mfp_2$, respectively.

Suppose that $\ga_2=\ga_1$ or $\ga_2=1-\ga_1$. If $M$ admits an Einstein adapted metric $g_M$ such that $g_F$ is not
Einstein, then the corresponding bisymmetric triple
$(\mfg,\mfk,\mfl)$ is one of the triples in Table
\ref{nonbimII}. The Einstein adapted metrics are, up to homothety, given by,

$$g_{\mfm}=\frac{1}{X_1}B_{\mfp_1}\oplus\frac{1}{X_2}B_{\mfp_2}\oplus B_{\mfn}$$

and $X_1$, $X_2$ are indicated in Table \ref{nonbimII}.\ethm

\li

This paper is organized as follows: in Section \ref{sectionsymF} we state some results about Einstein adapted metrics on the total space of a homogeneous fibration with symmetric fiber. In Section \ref{cops} we provide formulae to compute the eigenvalues of Casimir operators intervenient in the Einstein equations of an adapted metric. We prove all the results stated above for bisymmetric fibrations in Section \ref{bsSection}. In Section \ref{4symsection}, we obtain a list of all compact simply-connected $4$-symmetric spaces of maximal rank which admit an Einstein adapted metric. The results presented in Section \ref{4symsection} are a immediate consequence of the results for bisymmetric fibrations proved in this paper and the classification of $4$-symmetric spaces obtained by Jimenez in \cite{Ji}. Finally, all the tables mentioned in this paper are presented in Section \ref{alltabs}.

\section{Riemannian Fibrations with Symmetric
Fiber}\label{sectionsymF}

Let $G$ be a compact connected semisimple Lie group and $L\varsubsetneq K\varsubsetneq G$
connected closed non-trivial subgroups such that $N=G/K$ is isotropy irreducible and $F= K/L$ is a simply-connected symmetric space. We consider the natural fibration

\beq \label{symfibration} M=G/L\rightarrow G/K=N\textrm{ with fiber }F=K/L\eeq

and investigate the existence of Einstein adapted metrics on $M$. We use some results proved in \cite{Fa2} to deduce conditions for existence of an Einstein adapted metric and to describe some special class of Einstein adapted metrics on $M$, in the case of a fibration given by (\ref{symfibration}). We first introduce some notation and recall some useful results from \cite{Fa2}.

Let $B=-\Phi$, where $\Phi$ is the Killing form of $G$. We recall that since $G$ is compact and semisimple, $B$ is positive definite. We consider a $B$-orthogonal decomposition

\beq\label{decFsym} \mfg = \mfl\oplus\mfm = \mfl\oplus\mfp\oplus\mfn\eeq

of $\mfg$. We choose (\ref{decFsym}) such that $\mfg = \mfl\oplus\mfm$, $\mfg = \mfk\oplus\mfn$ and $\mfk=\mfl\oplus\mfp$ are reductive decompositions for $M$, $N$ and $F$, respectively. Since $F=K/L$ is a compact simply-connected symmetric space, we consider its DeRham
decomposition

\beq \label{deRham1}K/L=K_1/L_1\times\ldots\times K_s/L_s,\eeq

where each $K_a$ is simple. In particular, each $K_a/L_a$ is an irreducible symmetric space. By $\mfk_a$ and $\mfl_a$ we denote the Lie algebras of
$K_a$ and $L_a$, respectively. Let $\mfp_a$ be a symmetric reductive complement of $\mfl_a$ in
$\mfk_a$. Hence, we have a decomposition of $\mfm$ given by

\beq \mfm=\mfp_1\oplus\ldots\oplus\mfp_s\oplus\mfn,\eeq

where each $\mfp_a$ is $Ad\,L$-irreducible and $\mfn$ is $Ad\,K$-irreducible. Throughout, we assume that hypothesis (\ref{submodulesHyp}) is satisfied, i.e.,  $\mfp_1,\,\ldots\mfp_s$ are inequivalent $Ad\,L$-submodules and $\mfp,\,\mfn$ do not contain equivalent $Ad\,L$-submodules.

Under the construction above, any adapted metric $g_M$ on $M$ (see Section \ref{intro}) is defined by an $Ad\,L$-invariant Euclidean product on $\mfm$ of the form

\beq\label{mdef}g_{\mfm}=\left(\oplus_{a=1}^s\la_aB_{\mfp_a}\right)\oplus\mu B_{\mfn},\,\la_a,\,\mu>0.\eeq

Recall the Casimir operator of a subspace defined by (\ref{casdef}). Since each $\mfk_a$ is simple and $\mfn$ is an irreducible $Ad\,K$-module, the Casimir operator of $\mfk$
is scalar on each $\mfk_a$ and on $\mfn$. Let $\ga_a$ and $c_{\mfk,\mfn}$ be the corresponding eigenvalues:

\beqar\label{gasdef}C_{\mfk}\mid_{\mfk_a}=\ga_aId\\
 C_{\mfk}\mid_{\mfn}=c_{\mfk,\mfn}Id\eeqar

Necessary conditions for existence of Einstein adapted metrics are given as algebraic conditions on the Casimir operators of $\mfk$ and $\mfp_a$. More exactly, we recall the following result proved in \cite{Fa2}:

\bthm \label{cond3} \cite{Fa2} Let $M=G/L\rightarrow G/K=N$ be a homogeneous fibration, for a compact connected semisimple Lie group $G$ such that $N$ is isotropy irreducible. If $M$ admits an Einstein adapted
metric, then there are constants $\la_1,\ldots,\la_s>0$
such that $\sum_{a=1}^s\la_aC_{\mfp_a}$ is scalar on $\mfn$, where $C_{\mfp_a}$ is the Casimir operator of $\mfp_a$.
\ethm

The condition given in Theorem \ref{cond3} is a very useful test for existence of Einstein adapted metrics as we shall see in Section \ref{bsSection}. In particular, the condition in Theorem \ref{cond3} is satisfied if the Casimir operators $C_{\mfp_a}$ are scalar on $\mfn$. Under this assumption, Einstein adapted metrics on $M$ are given by positive solutions of a system of $s$ algebraic equations of degree $s-1$ with $s$ variables. More precisely we can state the following:

\bthm  \label{eqFsym} \cite[\S2.5]{Fa1} Let $M=G/L\rightarrow G/K$ be a homogeneous fibration with symmetric fiber $F=K/L$ as in (\ref{symfibration}). Suppose that $C_{\mfp_a}$ is scalar on $\mfn$, i.e., $C_{\mfp_a}\mid_{\mfn}=b_aId$, $a=1,\ldots,s$.

$M$ admits an Einstein adapted metric if and only if there are positive solutions of the
following system of $s$ algebraic equations on the unknowns
$X_1,\ldots,X_s$:\footnote{$\widehat{X_a}$ means that $X_a$ does
not occur in the product. }

\beqar 2\ga_1X_1^2X_a+(1-\ga_1)X_a-2\ga_aX_1X_a^2-(1-\ga_a)X_1=0,
\,a=2,\ldots,s\nonumber\\
2\sum_{a=1}^sb_aX_1\ldots\widehat{X_a} \ldots X_s-4rX_1\ldots
X_s+2\ga_1X_1^2X_2\ldots X_s+(1-\ga_1)X_2\ldots X_s=0,\nonumber
\eeqar

where $\ga_a$ is the eigenvalue of $C_{\mfk}$ on $\mfp_a$,
$r=\frac{1}{2}\left(\frac{1}{2}+c_{\mfk,\mfn}\right)$ and
$c_{\mfk,\mfn}$ is the eigenvalue of $C_{\mfk}$ on $\mfn$. To each
$s$-tuple $(X_1,\ldots,X_s)$ corresponds, up to homothety, an adapted metric on $M$ given by

$$g_{\mfm}=\oplus_{a=1}^s\frac{1}{X_a}B_{\mfp_a}\oplus
B_{\mfn}.$$ \ethm

The proof of Theorem \ref{eqFsym} is out of the scope of this paper and can be found in \cite[\S2.5]{Fa1}. A special class of adapted metrics is the class of binormal metrics. A metric is said to be \textbf{binormal} if it is defined by an $Ad\,L$-invariant Euclidean product of the form

\beq\label{mdefbi}g_{\mfm}=\la B_{\mfp}\oplus \mu B_{\mfn}.\eeq

Einstein binormal metrics are studied in \cite{Fa2}. We state the following two results proved in \cite{Fa2} which we use to describe Einstein binormal metrics when the fiber is a symmetric space.

\bthm \cite{Fa2}\label{binormal1} Let $M=G/L\rightarrow G/K$ be a homogeneous fibration, for a compact connected semisimple Lie group $G$ and $L\varsubsetneq K\varsubsetneq G$ connected closed non-trivial subgroups of $G$. Let $\mfp=\mfp_1\oplus\ldots\oplus\mfp_s$ be a $B$-orthogonal decomposition into irreducible $Ad\,L$-submodules and $\mfn=\mfn_1\oplus\ldots\oplus\mfn_n$ be a $B$-orthogonal decomposition into irreducible $Ad\,K$-submodules.

\li

(i) If $C_{\mfp}$ is not scalar on some $\mfn_j$, then there are no Einstein binormal metrics on $M$;

(ii) Suppose that $C_{\mfp}$ is scalar on each $\mfn_j$, i.e., $C_{\mfp}\mid_{\mfn_j}=b^jId_{\mfn_j}$. Then
there is a one-to-one correspondence, up to homothety, between Einstein binormal metrics on $M$ and positive solutions of the
following set of quadratic equations on the unknown $X\in\reals$:

\beqar \label{einI1}\de_{ij}^{\mfk}(1-X)=\de_{ij}^{\mfl}, \textrm{ if
}n>1, \\
\label{einI2}(2\de_{ab}^{\mfl}+\de_{ab}^{\mfk})X^2=\de_{ab}^{\mfk},
\textrm{ if } s>1, \\
\label{einI3}\left(\ga_a+2c_{\mfl,a}\right)X^2-\left(1+2c_{\mfk,j}\right)X+(1-\ga_a+2b^j)=0,\eeqar

for $a,b=1,\ldots,s$ and $i,j=1,\ldots,n$, where $c_{\mfl,a}$
is the eigenvalue of $C_{\mfl}$ on $\mfp_a$, $\ga_a$ is the constant determined
by

$$\Phi_{\mfk}\mid_{\mfp_a\times\mfp_a}=\ga_a\Phi_{\mfp_a},$$

$c_{\mfk,j}$ is the eigenvalue of $C_{\mfk}$ on $\mfn_j$ and the
$\de$'s are the differences $\de_{ij}^{\mfk}=c_{\mfk,i}-c_{\mfk,j}$, $\de_{ij}^{\mfl}=c_{\mfl,i}-c_{\mfl,j}$, $\de_{ab}^{\mfk}=\ga_a-\ga_b$ and $\de_{ab}^{\mfl}=c_{\mfl,a}-c_{\mfl,b}$.

If such a positive solution $X$ exists, then Einstein binormal metrics are, up to
homothety, defined by $$g_{\mfm}=B_{\mfp}\oplus XB_{\mfn}.$$
\ethm

\bthm \cite{Fa2}\label{binormal4} Let $M=G/L\rightarrow G/K$ be a homogeneous fibration, for a compact connected semisimple Lie group $G$ and $L\varsubsetneq K\varsubsetneq G$ connected closed non-trivial subgroups of $G$. Suppose $F$ is not isotropy irreducible and
that there exists a constant $\al$ such that

$$\Phi\circ C_{\mfl}\mid_{\mfp\times\mfp}=\al
\Phi_{\mfk}\mid_{\mfp\times\mfp}.$$

For $a=1,\ldots,s$, let $\ga_a$ be the constant determined by
\beq\label{defgagen}\Phi_{\mfk}\mid_{\mfp_a\times\mfp_a}=\ga_a\Phi\mid_{\mfp_a\times\mfp_a}.\eeq

If for some $a,b=1,\ldots,s$, $\ga_a\neq \ga_b$ and there exists on
$M$ an Einstein binormal metric, then the number $\sqrt{2\al+1}$ is
a rational. \ethm

From Theorems \ref{binormal1} and \ref{binormal4} we deduce the following result for the case when the fiber is a symmetric space and the base is isotropy irreducible:

\bcor \label{binormal5} Let $M=G/L\rightarrow G/K$ be a homogeneous fibration with symmetric fiber $F=K/L$ as in (\ref{symfibration}).

(i) If $C_{\mfp}$ is not scalar on $\mfn$ or $C_{\mfk}$ is not
scalar on $\mfp$, then there is no Einstein binormal metric on
$M$.

(ii) Suppose that $C_{\mfp}$ is scalar on $\mfn$ and $C_{\mfk}$ is
scalar on $\mfp$, i.e., $C_{\mfp}\mid_{\mfn}=bId$ and
$C_{\mfk}\mid_{\mfp}=\ga Id$. There is an one-to-one
correspondence between Einstein binormal  metrics on $M$ and
positive roots of the quadratic equation

\beq\label{polbin5}2\ga
X^2-\left(1+2c_{\mfk,\mfn}\right)X+(1-\ga+2b)=0.\eeq

where $c_{\mfk,\mfn}$ is the eigenvalue of $C_{\mfk}$ on $\mfn$. If such
a positive solution $X$ exists, then Einstein binormal metrics
are, up to homothety, given by

$$g_{\mfm}=B_{\mfp}\oplus XB_{\mfn}.$$
\ecor

\bproof We first note that the constant $\ga_a$ defined by (\ref{defgagen}) coincides with the eigenvalue of $\mfk$ on $\mfk_a$ as defined in (\ref{gasdef}). If $F$ is isotropy irreducible, then $C_{\mfk}$ is scalar on $\mfp$. Suppose $F$ is not irreducible. Since $F$ is a symmetric space, then
$\Phi\circ C_{\mfl}\mid_{\mfp\times\mfp}=\al
\Phi_{\mfk}\mid_{\mfp\times\mfp}$, with $\al=\frac{1}{2}$. The
number $\sqrt{2\al+1}=\sqrt{2}$ is not a rational. Hence, from Theorem \ref{binormal4} we conclude that if exists a binormal
Einstein metric on $M$, then $\ga_1=\ldots=\ga_s=\ga$ for some
constant $\ga$, i.e., the Casimir operator of $\mfk$ is scalar on $\mfp$.

By using Theorem \ref{binormal1}, the condition that $C_{\mfp}$ is scalar on $\mfn$ is a necessary condition
for the existence of an Einstein binormal metric. Also, the condition (\ref{einI2}) from Theorem \ref{binormal1} is
satisfied since for $\ga_1=\ldots=\ga_s$, we have
$\de_{ab}^{\mfk}=\de_{ab}^{\mfl}=0$. Since $N=G/K$ is isotropy irreducible, (\ref{einI1}) is trivial. Finally, the polynomial (\ref{polbin5}) is just (\ref{einI3}) from Theorem \ref{binormal1}, for $c_{\mfl,a}=\frac{\ga_a}{2}=\frac{\ga}{2}$ and $n=1$.

$\Box$\eproof

If both fiber and base are isotropy irreducible symmetric spaces we obtain the following simplification of Corollary \ref{binormal5}.

\bcor \label{binSym1} Let $M=G/L\rightarrow G/K=N$ be a homogeneous fibration with symmetric fiber $F=K/L$ as in (\ref{symfibration}). Suppose that $F$ and $N$ are irreducible
symmetric spaces and $dim\,F>1$. There exists on $M$ an Einstein
adapted metric if and only if $C_{\mfp}$ is scalar on $\mfn$ and
$\triangle \geq 0$, where
$$\triangle=1-2\ga(1-\ga+2b),$$

$\ga$ is the eigenvalue of $C_{\mfk}$ on $\mfp$ and $b$ is the
eigenvalue of $C_{\mfp}$ on $\mfn$. If these two conditions are
satisfied, then Einstein adapted metrics are, up to homothety, determined by

$$g_{\mfm}=B_{\mfp}\oplus XB_{\mfn}\textrm{ where } X=\frac{1\pm\sqrt{\triangle}}{2\ga}.$$\ecor

\bproof It follows immediately from Corollary \ref{binormal5} and from the fact
that if $N$ is a symmetric space, then $c_{\mfk,\mfn}=\dfrac{1}{2}$.

$\Box$\eproof

Another special class of Einstein adapted metrics are those whose restriction $g_F$ to the fiber $F$ is also an Einstein metric. We state the following result from \cite{Fa2}:

\bthm \label{bfein} Let $g_M$ be an Einstein adapted metric on the homogeneous fibration $M=G/L\rightarrow G/K=N$, such that $N$ is isotropy irreducible, defined by the $Ad\,L$-invariant Euclidean product $g_{\mfm}=\left(\oplus_{a=1}^s\la_aB_{\mfp_a}\right)\oplus\mu B_{\mfn}.$

If $g_F$ is also Einstein, then

\beq
\frac{\la_a}{\la_b}= \frac{C_{\mfn,b}}{C_{\mfn,a}}\eeq

for $,\,a,b=1,\ldots,s$, where $c_{\mfn,a}$ is defined by $\Phi(C_{\mfn}\cdot,\cdot)\mid_{\mfp_a\times\mfp_a}=c_{\mfn,a}\Phi\mid_{\mfp_a\times\mfp_a}$.

In particular, there exists at most one $K$-invariant metric $g_F$ on $F$ such that $g_M$ is Einstein.

\ethm

Below we obtain a necessary condition for existence of an Einstein adapted metric with this property in terms of the Casimir operator of $\mfk$, when the fiber is a symmetric space.

\bcor \label{gfeincor} Let $M=G/L\rightarrow G/K=N$ be a homogeneous fibration with symmetric fiber $F=K/L$ as in (\ref{symfibration}).  Let $\ga_a$ be the eigenvalue of $C_{\mfk}$
on $\mfp_a$ and $g_M$ an adapted metric on $M$. If
$g_M$ and $g_F$ are both Einstein, then

$$\ga_a=\ga_b\textrm{ or }\ga_a=1-\ga_b,\,a,b=1,\ldots,s.$$
\ecor

\bproof If $F$ is irreducible, then the statement is trivial. We suppose that $F$ is reducible. First we recall that since $F$ is a symmetric space, the Ricci curvature of $g_F$ is given by $Ric^F=\frac{1}{2}\Phi_{\mfk}$, where $\Phi_{\mfk}$ is the Killing form of $\mfk$. Therefore, for $X\in\mfp_a$,

\beq \label{ricgF}Ric^F(X,X)=-\frac{1}{2}\Phi_{\mfk}(X,X)=-\frac{1}{2}\Phi(C_{\mfk}X,X)=-\frac{\ga_a}{2}\Phi(X,X).\eeq

If $g_F$ is Einstein with Einstein constant $E_F$, then (\ref{ricgF}) implies that $\frac{\ga_a}{2}=E\la_a$. Consequently, we obtain the following relations

\beq \label{rellaga}\frac{\la_a}{\la_b}=\frac{\ga_a}{\ga_b}\eeq

Let $C_{\mfn,a}$ be the constant defined by $\Phi(C_{\mfn}\cdot,\cdot)\mid_{\mfp_a\times\mfp_a}=c_{\mfn,a}\Phi\mid_{\mfp_a\times\mfp_a}$. Clearly, we have

\beq \label{relcnga} C_{\mfn,a}=1-\ga_a. \eeq

The identity given by Theorem \ref{bfein}, together with (\ref{rellaga}) and (\ref{relcnga}), imply that

$$\frac{\ga_a}{\ga_b}=\frac{1-\ga_b}{1-\ga_a}.$$

Consequently, $\ga_a=\ga_b$ or $\ga_a=1-\ga_b$.

$\Box$\eproof

For a homogeneous fibration with symmetric fiber, Einstein binormal metrics restrict to an Einstein metric on the fiber:

\bcor \label{binormal6} Let $M=G/L\rightarrow G/K$ be a homogeneous fibration with symmetric fiber $F=K/L$ as in (\ref{symfibration}). If
there exists on $M$ an Einstein binormal metric $g_M$, then $g_F$
is Einstein. The converse holds if $C_{\mfk}$ is scalar on
$\mfp$.\ecor

\bproof If $F$ is irreducible, then the statement is trivial. We suppose that $F$ is reducible. If $M$ admits an Einstein binormal metric $g_M$, then, by Corollary \ref{binormal5}, $C_{\mfk}$ is scalar on $\mfp$. Hence, $\ga_1=\ldots=\ga_s=\ga$ for some $\ga$. By (\ref{ricgF}) we obtain

$$Ric^F=\frac{\ga}{2}B_{\mfp}$$

and, consequently, $g_F$ is Einstein with Einstein constant $E_F=\frac{\ga}{2\la}.$

Conversely, let $g_M$ be any Einstein
adapted metric on $M$ as in (\ref{mdef}). If $C_{\mfk}$ is scalar on $\mfp$, then
$\ga_1=\ldots=\ga_s$. Hence, if $g_F$ is also Einstein, the equality (\ref{rellaga}) in the proof above implies that $\la_1=\ldots=\la_s$. Therefore, $g_M$ is binormal.

$\Box$\eproof

There might be non-binormal Einstein adapted metrics whose restriction to the fiber is still Einstein. The existence of an Einstein adapted metric $g_M$ such that $g_F$ is Einstein implies that the Casimir operator of $\mfk$ satisfies one of the condition given in Corollary \ref{gfeincor}. If in addition we suppose that $F$ is the direct product of only two isotropy irreducible symmetric spaces, the Einstein equations in Theorem \ref{eqFsym} can be solved. The result below describes Einstein adapted metrics such that $g_F$ is also Einstein in this case.

\bcor \label{gfeincor2} Let $M=G/L\rightarrow G/K$ be a homogeneous fibration with symmetric fiber $F=K/L$ as in (\ref{symfibration}). Suppose that $\mfp=\mfp_1\oplus\mfp_2$ is a decomposition of $\mfp$ into $Ad\,L$-irreducible submodules and $C_{\mfp_a}\mid_{\mfn}=b_aId$, for some constants $b_a$, $a=1,2$.  Let $\ga_a$ be the eigenvalue of $C_{\mfk}$ on $\mfp_a$ and $c_{\mfk,\mfn}$ be the eigenvalue of $C_{\mfk}$ on
$\mfn$.

If there exists on $M$ an Einstein adapted metric $g_M$ such that
$g_F$ is also Einstein, then one of the following cases holds:

(i) $\ga_2=\ga_1$ and $\triangle\geq 0$, where

$$\triangle=(1+c_{\mfk,\mfn})^2-8\ga_1(1-\ga_1+2b).$$

If these two conditions are satisfied, then $g_M$ is
the binormal metric given, up to homothety, by

$$g_{\mfm}=B_{\mfp}\oplus
XB_{\mfn},\,
where\,X=\frac{1+c_{\mfk,\mfn}\pm\sqrt{\triangle}}{2\ga_1}.$$

(ii) $\ga_2=1-\ga_1$ and $D(\ga_1)\geq 0$, where

$$D(\ga_1)=4r^2-4b_1\ga_1-4b_2(1-\ga_1)-2\ga_1(1-\ga_1)$$

and $r=\frac{1}{2}\left(\frac{1}{2}+c_{\mfk,\mfn}\right)$. If
these two conditions are satisfied, then $g_M$ is given, up
to homothety, by

$$g_{\mfm}=\frac{1}{X_1}B_{\mfp_1}\oplus\frac{1}{X_2}B_{\mfp_2}\oplus
B_{\mfn},$$

where $$X_2=\frac{\ga_1 X_1}{1-\ga_1}\textrm{  and
}X_1=\frac{2r\pm\sqrt{D(\ga_1)}}{2\ga_1}.$$\ecor

\bproof Let $g_M$ be an Einstein adapted metric on $M$ associated to the $Ad\,L$-invariant Euclidean product $g_{\mfm}=\left(\oplus_{a=1}^s\la_aB_{\mfp_a}\right)\oplus\mu B_{\mfn}$.

If the restriction to $F$ is also Einstein then, by Corollary \ref{gfeincor}, either $\ga_2=\ga_1$ or $\ga_2=1-\ga_1$. The statement (i), for the case $\ga_2=\ga_1$, follows from Corollaries
\ref{binormal5} and \ref{binormal6}. In the case $\ga_2=1-\ga_1$, the equality (\ref{rellaga}) implies that

\beq\label{rellambdas}\frac{\la_1}{\la_2}=\frac{\ga_1}{1-\ga_1}.\eeq

On the other hand, by simplifying the equations given by Theorem \ref{eqFsym}, an adapted
Einstein metric on $M$ corresponds to positive solutions of the
equations

\beqar
\label{eqII11}2\ga_1X_1^2X_2+(1-\ga_1)X_2-2\ga_2X_1X_2^2-(1-\ga_2)X_1=0
\\
\label{eqII12}2b_1X_2+2b_2X_1-4rX_1X_2+2\ga_1X_1^2X_2+(1-\ga_1)X_2=0
\eeqar

where $X_a=\frac{\mu}{\la_a}$, for $a=1,2$. By using the identity
(\ref{rellambdas}), we write

\beq\label{x2} X_2=
\frac{\ga_1}{1-\ga_1}X_1.\eeq

We solve the system of equations
(\ref{eqII11}) and (\ref{eqII12}) using (\ref{x2}) and $\ga_2=1-\ga_1$. This proves statement (ii).

$\Box$ \eproof

Even under the assumption that the fiber has only two isotropy subspaces, to classify all the Einstein adapted metrics is a complicated problem. The Einstein equations given by Theorem \ref{eqFsym} are still unmanageable. Under the conditions $\ga_2=\ga_1$ or $\ga_2=1-\ga_1$, on the Casimir operator of $\mfk$, it is possible to solve these equations. Corollaries \ref{contcor1} and \ref{contcor2} below describe all the Einstein adapted metrics on these spaces. Their proofs are similar to the proof of Corollary \ref{gfeincor2} by solving the equations given by Theorem \ref{eqFsym}. We omit these two proofs which can be found in \cite[\S 2.5]{Fa1}. We note that, in particular, if $M$ admits an Einstein adapted metric whose restriction to the fiber is Einstein, then the following two results give all the others Einstein adapted metrics on $M$.

\bcor \cite[p.50]{Fa1} \label{contcor1} Let $M=G/L\rightarrow G/K$ be a homogeneous fibration with symmetric fiber $F=K/L$ as in (\ref{symfibration}). Suppose that $\mfp=\mfp_1\oplus\mfp_2$ is a decomposition of $\mfp$ into $Ad\,L$-irreducible submodules and $C_{\mfp_a}\mid_{\mfn}=b_aId$, for some constants $b_a$, $a=1,2$.  Let $\ga_a$ be the eigenvalue of $C_{\mfk}$ on $\mfp_a$ and $c_{\mfk,\mfn}$ be the eigenvalue of $C_{\mfk}$ on
$\mfn$.

Suppose that $\ga_2=\ga_1$, i.e., $C_{\mfk}$ is scalar on $\mfp$.
If there exists on $M$ an Einstein adapted metric $g_M$, then one
of the following two cases holds:

(i) $g_F$ is also Einstein and $g_M$ is a binormal metric given by
Corollary \ref{gfeincor2} (i).

(ii) $D(\ga_1)\geq 0$, where

$$D(\ga_1)=4r^2(1-\ga_1)-2\ga_1(2b_2+1-\ga_1)(2b_1+1-\ga_1)$$

and $r=\frac{1}{2}\left(\frac{1}{2}+c_{\mfk,\mfn}\right)$. The
metric $g_M$ is given, up to homothety, by

$$g_{\mfm}=\frac{1}{X_1}B_{\mfp_1}\oplus\frac{1}{X_2}B_{\mfp_2}\oplus
B_{\mfn},$$

where $$X_2=\frac{1-\ga_1}{2\ga_1 X_1}\textrm{  and
}X_1=\frac{2r(1-\ga_1)\pm\sqrt{(1-\ga_1)D(\ga_1)}}{2\ga_1(2b_2+1-\ga_1)}.$$
In this second case, $g_F$ is not Einstein and $g_M$ is not
binormal.\ecor

\bcor \cite[p.50]{Fa1} \label{contcor2} Let $M=G/L\rightarrow G/K$ be a homogeneous fibration with symmetric fiber $F=K/L$ as in (\ref{symfibration}). Suppose that $\mfp=\mfp_1\oplus\mfp_2$ is a decomposition of $\mfp$ into $Ad\,L$-irreducible submodules and $C_{\mfp_a}\mid_{\mfn}=b_aId$, for some constants $b_a$, $a=1,2$.  Let $\ga_a$ be the eigenvalue of $C_{\mfk}$ on $\mfp_a$ and $c_{\mfk,\mfn}$ be the eigenvalue of $C_{\mfk}$ on
$\mfn$.

Suppose that $\ga_2=1-\ga_1$. If there exists on $M$ an Einstein
adapted metric $g_M$, then one of the following two cases holds:

(i) $g_F$ is also Einstein and $g_M$ is the metric given by
Corollary \ref{gfeincor2} (ii).

(ii) $D(\ga_1)\geq 0$, where

$$D(\ga_1)=4r^2-2(2b_2+\ga_1)(2b_1+1-\ga_1)$$

and $r=\frac{1}{2}\left(\frac{1}{2}+c_{\mfk,\mfn}\right)$. The
metric $g_M$ is given, up to homothety, by

$$g_{\mfm}=\frac{1}{X_1}B_{\mfp_1}\oplus\frac{1}{X_2}B_{\mfp_2}\oplus
B_{\mfn},$$

where
$$X_2=\frac{1}{2X_1}\textrm{  and
}X_1=\frac{2r\pm\sqrt{D(\ga_1)}}{2(2b_2+\ga_1)}.$$

$g_M$ is never binormal and in the second case $g_F$ is not
Einstein.\ecor

\section{The Casimir Operators}\label{cops}

In this section we consider a compact connected semisimple Lie group $G$ and $L\varsubsetneq K\varsubsetneq G$ connected closed non-trivial
subgroups. We only suppose that $L$ and $K$ have maximal rank in $G$ and $N=G/K$ is isotropy irreducible. Using the notation from previous sections, we provide formulae to compute the eigenvalues of the Casimir operator of the vertical space $\mfp$ on the horizontal direction $\mfn$ and the eigenvalues of the Casimir operator of $\mfk$ on the vertical direction $\mfp$. Also as above,

$$\mfp=\mfp_1\oplus\ldots\oplus\mfp_s$$

is a decomposition of $\mfp$ into $Ad\,L$-irreducible submodules. We recall some theory of roots of a semisimple Lie group (see \cite{He},\cite{Hu}). Let $\mfh$ be a Cartan subalgebra of $\mfg^{\complex}$,
such that $\mfh\subset\mfl^{\complex}$, and let $\R$ be a system of
nonzero roots for $\mfg^{\complex}$, with respect to $\mfh$. We have a decomposition of $\mfg^{\complex}$ into root subspaces given by

\beq \mfg^{\complex}=\mfh\oplus(\oplus_{\al\in \R}\mfg^{\al}).\eeq

We consider a standard
normalized basis $\{E_{\al}\}_{\al\in \R}$ of $\oplus_{\al\in R}\mfg^{\al}$ and the elements $H_{\al}=[E_{\al},E_{-\al}]\in \mfh$. We recall that $\Phi(H_{\al},h)=\al(h)$, for
every $h\in\mfh$. In particular, the length $|\al|$ of a root
$\al\in\R$ in $\mfg$ is defined by $ |\al|^2=\al(H_{\al})=\Phi(H_{\al},H_{\al})$. For every $\al,\be\in\R$ such that $\al+\be\in\R$, the \emph{structure constants} are the numbers
$N_{\al,\be}\in \complex$ defined by

\beq\label{sc1}[E_{\al},E_{\be} ]=N_{\al,\be}E_{\al+\be},\eeq

They satisfy the properties

\beq \label{scprop}N_{\al,\be}=-N_{\be,\al}\textrm{ and }N_{-\al,\be+\al}=N_{-\be,-\al}=N_{\al,\be}.\eeq

We identify $\mfg$ with the compact real form generated the by maximal toral subalgebra $i\mfh_{\reals}$ and by the elements
\beq\label{basiscf}
X_{\al}=\frac{E_{\al}-E_{-\al}}{\sqrt{2}}\textrm{ and
}Y_{\al}=\frac{i(E_{\al}+E_{-\al})}{\sqrt{2}}, \,\al\in\R^+,\eeq
where $\R^+\subseteq \R$ is a subset of positive roots. We define the subsets of roots

\beq \R_{\mfk}=\{\al\in\R:E_{\al}\in\mfk^{\complex}\}\textrm{ and }\R_{\mfp_a}=\{\al\in\R:E_{\al}\in\mfp_a^{\complex}\}.\eeq

Since $\mfl$ has maximal rank and $\mfk=\mfl\oplus\mfp$,

\beq \mfp_a^{\complex}=<E_{\al},E_{-\al} :\al\in\R_{\mfp_a}^+>.\eeq

Moreover, since $\Phi(E_{\al},E_{-\al})=1$, the bases
$\{E_{\al}\}_{\al\in\R_{\mfp_a}}$ and
$\{E_{-\al}\}_{\al\in\R_{\mfp_a}}$ of $\mfp_a^{\complex}$ are dual
with respect to $\Phi$. Consequently, the Casimir
operator of $\mfp_a^{\complex}$ may be written as

\beq
\label{cpacomplex}C_{\mfp_a^{\complex}}=\sum_{\al\in\R_{\mfp_a}}ad_{E_{\al}}ad_{E_{-\al}}.\eeq

Also, since $\mfk$ has maximal rank and $\mfg=\mfk\oplus\mfn$, we have
$\mfn^{\complex}=<E_{\al}:\al\in\R_{\mfn}>$, where

\beqar \R_{\mfn}=\{\al\in\R:E_{\al}\in \mfn^{\complex}\}=\R-\R_{\mfk}.\eeqar

By hypothesis the $Ad\,K$-module $\mfn$ is irreducible. If
$\mfn=\oplus_j\mfn^j$ is a decomposition of $\mfn$ into
irreducible $Ad\,L$-submodules, we write

\beq \R_{\mfn^j}=\{\phi\in\R: E_{\phi}\in(\mfn^j)^{\complex}\}.\eeq

\li

We recall that a necessary condition for existence of an Einstein
adapted metric  on $M$, given by Theorem \ref{cond3}, is that
there are $\la_1,\ldots,\la_s>0$ such that the operator
$\sum_{a=1}^s\la_aC_{\mfp_a}$ is scalar on $\mfn$. The Casimir
operator $C_{\mfp_a}$ is necessarily scalar on the irreducible
$Ad\,L$-submodules $\mfn^j$. We
shall write $b_a^j$ for this eigenvalue, i.e.,

\beq\label{cpa1}
C_{\mfp_a}\mid_{\mfn^j}=b_a^j Id.\eeq

Furthermore, the eigenvalue of $C_{\mfp_a}$ on $\mfn^j$,
$b_a^j$, must coincide with the eigenvalue of
$C_{\mfp_a^{\complex}}$ on $(\mfn^j)^{\complex}$. Hence, we also
have

\beq \label{cpa2}
b_a^j=\Phi(C_{\mfp_a^{\complex}}E_{\phi},E_{-\phi}),\,\forall\,\phi\in\R_{\mfn^j}.\eeq

We next prove a formula to compute the eigenvalues $b_a^j$. For any roots $\phi$ and $\al$ let $\phi+n\al$, $p_{\al\phi}\leq
n\leq q_{\al\phi}$, be the $\al$-series containing $\phi$ (see \cite[Chap.III \S
4]{He}). By definition, the $\al$-series containing $\phi$ is the set of all
roots of the form $\phi+n\al$ where $n$ is an integer. For $\al,\,\phi\in\R$,

\beq\label{sc1}N_{\al,\phi}^2=\frac{q_{\al\phi}(1-p_{\al\phi})}{2}\al(H_{\al}).\eeq

\bprop \label{bejs} Let $\mfn^j\subseteq\mfn$ be an irreducible $Ad\,L$-submodule and $\phi\in \R_{\mfn^j}$. For $a=1,\ldots,s$,

$$b_a^j=\dfrac{1}{2}\sum _{\al\in
\R_{\mfp_a}^+}d_{\al\phi}|\al|^2,$$

where $d_{\al\phi}=q_{\al\phi}-p_{\al\phi}-2p_{\al\phi}q_{\al\phi}$
and $\phi+n\al$, $p_{\al\phi}\leq n\leq q_{\al\phi}$ is the
$\al$-series containing $\phi$.

\eprop

\bproof By using (\ref{cpa2}) and (\ref{cpacomplex}) we obtain the
following:

$$\bar{rl}b_a^j= & \Phi(C_{\mfp_a^{\complex}}E_{\phi},E_{-\phi})\\ \\

= &
\sum_{\al\in\R_{\mfp_a}}\Phi([E_{-\al},[E_{\al},E_{\phi}]],E_{-\phi})\\
\\

= &
\sum_{\al\in\R_{\mfp_a}}N_{\al,\phi}\Phi([E_{-\al},E_{\phi+\al}],E_{-\phi})\\
\\
= &
\sum_{\al\in\R_{\mfp_a}}N_{\al,\phi}N_{-\al,\phi+\al}\Phi(E_{\phi},E_{-\phi})\\
\\

= & \sum_{\al\in\R_{\mfp_a}}N_{\al,\phi}N_{-\al,\phi+\al}\ear$$

From (\ref{scprop}) we have $N_{-\al,\phi+\al}=N_{\al,\phi}$ and
we get

\beq \label{prel1}b_a^j=\sum_{\al\in\R_{\mfp_a}}N_{\al,\phi}^2=\sum_{\al\in\R_{\mfp_a}^+}(N_{\al,\phi}^2+N_{-\al,\phi}^2).\eeq

Now let $\phi+n\al$, $p_{\al\phi}\leq n\leq q_{\al\phi}$, be the
$\al$-series containing $\phi$. The number $N_{\al,\phi}^2$ is given by (\ref{sc1}). To compute $N_{-\al,\phi}^2$ we need the
$(-\al)$-series containing $\phi$. Clearly, this series is
$\phi-n'\al$, where $-q_{\al\phi}\leq n'\leq -p_{\al\phi}$. Hence,

\beq \label{prel2}N_{-\al,\phi}^2=\frac{-p_{\al\phi}(1-(-q_{\al\phi}))}{2}(-\al)(H_{-\al})=\frac{-p_{\al\phi}(1+q_{\al\phi})}{2}\al(H_{\al}).\eeq

Replacing (\ref{sc1}) and (\ref{prel2}) in (\ref{prel1}) concludes the proof.

$\Box$\eproof

\li

Let us consider a decomposition of $\mfk$ into its center $\mfk_0$
and simple ideals $\mfk_a$, for $a=1,\ldots, t$,

\beq\label{dec1}\mfk=\mfk_0\oplus\mfk_1\oplus\ldots\oplus\mfk_t,\eeq

and let $\ga_a$ denote the eigenvalue of the Casimir operator of
$\mfk$ on $\mfk_a$. We present a formula to compute the eigenvalues
$\ga_a$'s by making use of dual Coxeter numbers (see \cite[\S V.5]{BD},\cite[10.4]{Hu}). The \textbf{dual Coxeter number} of a simple Lie algebra $\mfg$ is
the number given by

\beq h^*(\mfg)=\dfrac{1}{|\al|^2},\eeq

where $\al$ is a long root (see \cite{Pa}). We may suppose that $\mfh_a=\mfh\cap\mfk_a$ is a Cartan subalgebra
of $\mfk_a$ and thus a root of $\mfk_a$ can be viewed as a root for
$\mfg$. Hence we can compare lengths of roots of $\mfg$ with lengths
of roots of $\mfk_a$. So let $\de_a$ be the ratio of the square
length of a long root for $\mfg$ to that of $\mfk_a$, i.e.,

\beq \de_a=\dfrac{|\al|^2_{\mfg}}{|\be|^2_{\mfg}}=\dfrac{\Phi(H_{\al},H_{\al})}{\Phi(H_{\be},H_{\be})},\eeq

where $\al$ is a long root of $\mfg$ and $\be$ is a long root of
$\mfk_a$. Clearly, $\de_a=1$ if there exists only one length for
$\mfg$ or if both $\mfg$ and $\mfk_a$ have two lengths. If
$\de_a\neq 1$, $\de_a$ is equal to
either 2, if $G$ is of type $G_2$, or 3, if $G$ is of type $F_4$, $B_n$ or $C_n$. We state the following result by D. Panyushev which gives a formula to compute the eigenvalue of $C_{\mfk}$ on $\mfk_a$:

\bprop \label{gajs}\cite{Pa} Suppose that $\mfg$ is simple. Then

$$\ga_a=\dfrac{h^*(\mfk_a)}{\de_a . h^*(\mfg)},\,a=1,\ldots,s,$$

where $h^*(\mfk_a)$ and $h^*(\mfg)$ are the dual Coxeter numbers of
$\mfk_a$ and $\mfg$, respectively. \eprop

\section{Bisymmetric Fibrations}\label{bsSection}

In this section we prove Lemma \ref{bsclass} of classification of bisymmetric triples and show all the other results stated in Section \ref{intro} about classification of Einstein adapted metrics for bisymmetric fibrations.

We recall that for a bisymmetric fibration $M=G/L\rightarrow G/K=N$, as introduced in Section \ref{intro} (see (\ref{bsfibration})), the fiber $F=K/L$ is a compact simply-connected symmetric space. Therefore, we apply the results from Section \ref{sectionsymF} to investigate the existence of Einstein adapted metrics on $M$. Furthermore, the eigenvalues $b_a^j$'s and
$\ga_a$'s as in Section \ref{sectionsymF} are calculated using Propositions \ref{bejs} and \ref{gajs}. Their
values are indicated in Tables \ref{eigIexc}, \ref{eigIclass},
\ref{eigIIexc} and \ref{eigIIclass} in Section \ref{alltabs} for
each irreducible bisymmetric triple of maximal rank. Whereas the computation of the $\ga_a$'s is a straightforward application of Proposition \ref{gajs} together with the Coxeter numbers given in Table \ref{tabcoxeter}, the computations of the $b_a^j$'s are long. These calculations are outlined in \cite[Appendix A]{Fa1}, where all the essential information is presented.

\subsection{Classification}\label{class}

We prove Lemma \ref{bsclass} stated in Section \ref{intro} which classifies irreducible bisymmetric triples of maximal rank into type I or II.

A classification of isotropy irreducible symmetric spaces can be found in \cite{He}. By using this we obtain a
list of all possible triples $(\mfg,\mfk,\mfl)$ such that $\mfg$ is simple, $L\varsubsetneq K$ are subgroups of maximal rank of $G$ and $K/L$, $G/K$ are symmetric spaces with $N=G/K$ isotropy irreducible and $F=K/L$ simply-connected.

\li

\textbf{\emph{Proof of Lemma \ref{bsclass}:}} By inspection
of the classification of symmetric pairs $(\mfg,\mfk)$ of compact
type in \cite{He} we obtain that those of maximal rank are the pairs
in Tables \ref{spexc} and \ref{spclass}. We obtain the list of all bisymmetric triples by combining all possible $K$'s and $L$'s which make $K/L$ a symmetric space as well. We just need to observe that $K/L$ is isotropy irreducible or has two isotropy irreducible submodules.

We observe that the cases when $\mfk$ is the centralizer of a torus
are only the cases

$(\mfe_6,\mfso_{10}\oplus \reals )$,
$(\mfe_7,\mfe_6\oplus\reals)$, $(\mfso_{2n},\mfu_n)$,
$(\mfso_n,\reals\oplus\mfso_{n-2})$, $(\mfsp_{n},\mfu_n)$ and
$(\mfsu_n,\mfsu_p\oplus\mfsu_{n-p}\oplus\reals)$. In all the other cases $\mfk$ is semisimple. If $\mfk$
is simple then $(\mfk,\mfl)$ shall be an irreducible symmetric
pair, i.e., $\mfp$ is an irreducible $Ad\,L$-submodule. Thus,
$(\mfg,\mfk,\mfl)$ is of type I. In the cases where
$\mfk=\mfk_1\oplus\reals$ with $\mfk_1$ a simple ideal of $\mfk$,
since we require $\mfl$ to be of maximal rank, we have
$\mfl=\mfl_1\oplus\reals$, where $\mfl_1$ is a subalgebra of
$\mfk_1$ with maximal rank and $(\mfk,\mfl)\iso(\mfk_1,\mfl_1)$ is
an irreducible symmetric pair. Thus, in this case,
$\mfp\iso\mfp_1$ is also an irreducible $L$-invariant subspace and
$(\mfg,\mfk,\mfl)$ is of type I. In the cases where
$\mfk=\mfk_1\oplus\mfk_2$, with $\mfk_1$ and $\mfk_2$ simple
ideals of $\mfk$, we have $\mfl=\mfl_1\oplus\mfl_2$, where, for
$i=1,2$, $\mfl_i$ is a subalgebra of $\mfk_i$ of maximal rank.
Clearly, one of the $\mfl_i$'s must be proper as we require that
$\mfl$ is a proper subalgebra of $\mfk$. If both $\mfl_1$ ad
$\mfl_2$ are proper, then $\mfp=\mfp_1\oplus\mfp_2$, where
$\mfp_1$ and $\mfp_2$ are nonzero irreducible $\mfl$-invariant
subspaces. Hence, in this case, $(\mfg,\mfk,\mfl)$ is of type II.
If exactly one of $\mfl_i$'s coincides with $\mfk_i$, then
$(\mfk,\mfl)\iso(\mfk_j,\mfl_j)$ and $\mfp=\mfp_j$, for that $j$
satisfying $\mfl_j\neq\mfk_j$, and once again $(\mfg,\mfk,\mfl)$
is of type I. Finally, we have the case of the spaces
$(\mfsu_n,\mfsu_p\oplus\mfsu_{n-p}\oplus\reals)$,
$p=1,\ldots,n-1$. Clearly, $\mfl$ must be of the form
$\mfl=\mfl_1\oplus\mfl_2\oplus\reals$, where $\mfl_1$ and $\mfl_2$
are maximal rank subalgebras of $\mfsu_p$ and $\mfsu_{n-p}$,
respectively. We obtain a triple of type $I$ if exactly one of the
$\mfl_i$'s is proper and a triple of type II if both $\mfl_1$ and
$\mfl_2$ are proper.

$\Box$

\subsection{Einstein Adapted Metrics for Bisymmetric Triples of Type I}

We prove Theorem \ref{mtypeI} stated in Section \ref{intro}. The bisymmetric triples of type I are given in Tables \ref{eigIexc} and \ref{eigIclass} (see Lemma \ref{bsclass}).

We recall that for Type
I, the vertical isotropy subspace $\mfp$ is an irreducible $Ad\,L$-submodule. Since $F=K/L$ and $N=G/K$ are irreducible symmetric spaces, Einstein adapted metrics are given by Corollary \ref{binSym1}. In particular, any Einstein adapted metric is binormal.

We recall that $b^j$ is the eigenvalue of $C_{\mfp}$ on the $Ad\,L$-irreducible submodule $\mfn^j\subseteq\mfn$ (see (\ref{cpa1})). For each bisymmetric triple, these eigenvalues are indicated in Tables \ref{eigIexc} and \ref{eigIclass} in the column $b^j$. A necessary condition for existence of an Einstein adapted metric, given by Corollary \ref{binSym1}, is that $C_{\mfp}$ is scalar on $\mfn$. Hence, the first
test for existence of an adapted Einstein metric is to
observe if there exists only one value $b^j$ in the
corresponding columns of Tables \ref{eigIexc}  and
\ref{eigIclass}. If $C_{\mfp}$ is scalar on $\mfn$, we denote its unique eigenvalue by $b$ as in Corollary \ref{binSym1}.

The eigenvalue of the Casimir operator of $\mfk$ on $\mfp$ is denoted by $\ga$ as in Corollary \ref{binSym1}. For each bisymmetric triple, these eigenvalues are also indicated in Tables \ref{eigIexc} and \ref{eigIclass} in the column $\ga$.

\li

\textbf{\emph{Proof of Theorem \ref{mtypeI}:}} By Corollary
\ref{binSym1}, the existence of an adapted Einstein metric implies
that the Casimir operator of $\mfp$ is scalar on $\mfn$. By
inspection we conclude from Tables \ref{eigIexc} and \ref{eigIclass} that the only spaces satisfying this condition are
those corresponding to the labels

$$\ref{cpbn3},\,\ref{cpdn5},\,\ref{cpcn5},\,\ref{cpf41},\,\ref{cpf42}\,\ref{cpg21},\,\ref{cpg22},\,\ref{cpe81},\,\ref{cpe83},\,\ref{cpe88},\,\ref{cpe71},\,\ref{cpe74},\,\ref{cpe62},\,\ref{cpe63}$$

\ref{cpan1} for $l=\frac{p}{2}$, $p$ even ($b=\frac{p}{4n}$),

\ref{cpbn2} for $s=\frac{n-p}{2}$, $n-p$ even ($b=\frac{n-p}{2(2n-1)}$),

\ref{cpdn2} for $l=\frac{p}{2}$, $p$ even ($b=\frac{p}{4(n-1)}$),

\ref{cpcn2} for $l=\frac{p}{2}$, $p$ even ($b=\frac{p}{8(n+1)}$),

\ref{cpe78} for $p=1$ ($b=\frac{1}{9}$),

\ref{cpe64} for $p=1$ ($b=\frac{1}{8}$).

\li

For each of these cases, we compute $\triangle=1-2\ga(1-\ga+2b)$ given by Corollary \ref{binSym1}, and the values
obtained are as follows

{\footnotesize{
$$\bar{ll|ll}  & \triangle &  & \triangle \\ \hline

\xstrut\ref{cpan1} & \left(\frac{n-p}{n}\right)^2>0 & \ref{cpg22} & \frac{11}{18}>0 \\

\xstrut\ref{cpbn2} & \frac{4p^2+8p-4n+5}{(2n-1)^2}>0,\,\forall
p=\lfloor \frac{\sqrt{4n-1}}{2}\rfloor,\ldots,n-1 & \ref{cpe81} & \frac{7p^2-56p+113}{225}>0 \\

\xstrut\ref{cpbn3} & \left(\frac{2p+1}{2n-1}\right)^2>0 & \ref{cpe83} & \frac{196}{225}>0\\

\xstrut\ref{cpdn2} & \frac{p^2-(2n+1)p+n^2+1}{(n-1)^2}>0 & \ref{cpe88} & -\frac{2}{25}<0\\

\xstrut\ref{cpdn5} & \left(\frac{p-n}{n-1}\right)^2>0 & \ref{cpe71} & \frac{64}{81}>0
\\

\xstrut\ref{cpcn2} &
\frac{3p^2+(3-4n)p+2(n^2+1)}{2(n+1)^2}>0 & \ref{cpe74} & \frac{164-60p+5p^2}{324}>0,\,p=2,4\\

\xstrut\ref{cpcn5} & \left(\frac{n-p}{n+1}\right)^2>0 & \ref{cpe78} & \frac{25}{81}>0\\

\xstrut\ref{cpf41} &  \frac{106-63p+7p^2}{162}> 0,\,\textrm{ iff
}p=1,\,7 & \ref{cpe62} & \frac{1}{9}>0,\, p=2; -\frac{1}{9},\,p=4\\

\xstrut\ref{cpf42} & \frac{49}{81}>0 & \ref{cpe63} & \frac{25}{36}>0\\

\xstrut\ref{cpg21} & \frac{1}{4}>0 & \ref{cpe64} & \frac{1}{4}>0\\
\ear$$

 }}

\li

Finally, we compute $X=\frac{1\pm\sqrt{\triangle}}{2\ga}$ (see Corollary \ref{binSym1}) for those
cases when $\triangle>0$. The values of $X$ are indicated in Table
\ref{mIexc}. For each case, both $g_F$ and $g_N$ are Einstein because both $F$ and $N$ are isotropy irreducible.

$\Box$

\li

\brem The triples \ref{cpe88} and \ref{cpe74}, for $p=6$, do not admit an Einstein adapted metric since $\triangle<0$. These are the fibrations $$M=\frac{E_8}{SU(8)\times SU(2)}\rightarrow\frac{E_8}{E_7\times SU(2)}$$

and

$$M=\frac{E_7}{SO(6)\times SO(6)\times SU(2)}\rightarrow\frac{E_7}{SO(12)\times SU(2)}.$$

In both cases, the isotropy representation of $M$ has only two irreducible $Ad\,L$-modules, i.e., $\mfn$ is also $Ad\,L$-irreducible. Therefore, any $G$-invariant metric is an adapted metric. Hence, from Theorem \ref{mtypeI}, we conclude that $M$ has no $G$-invariant Einstein metrics. This conclusion is also obtained by Dickinson and Kerr in \cite{DK}.
\erem

\li

\subsection{Einstein Adapted Metrics for Bisymmetric Triples of Type II}

We prove Theorems \ref{biII}, \ref{gfeinII}, \ref{genIIexc} and \ref{genIIall} stated in Section \ref{intro}. The bisymmetric triples of type II are given in Tables \ref{eigIIexc} and \ref{eigIIclass} (see Lemma \ref{bsclass}).

For a bisymmetric fibration of type II, the fiber $F=K/L$ is the direct product of two irreducible symmetric spaces. Hence, the vertical subspace $\mfp$ decomposes as $\mfp=\mfp_1+\mfp_2$, where $\mfp_i$, $i=1,2$, are
irreducible $Ad\,L$-modules. Moreover, $\mfn$ is an irreducible $Ad\,K$-module. According to Theorem \ref{cond3}, a necessary condition for existence of an Einstein adapted metric is that there exist $\la_1,\,\la_2>0$ such that the operator

\beq \label{condIIb}\la_1C_{\mfp_1}+\la_2C_{\mfp_2}\textrm{ is scalar on }\mfn.\eeq

The eigenvalues $b_1^j$, $b_2^j$ of $C_{\mfp_1}$, $C_{\mfp_2}$, respectively, on the $Ad\,L$-irreducible submodules $\mfn^j\subseteq\mfn$ (see (\ref{cpa1})) are given in Tables \ref{eigIIexc} and \ref{eigIIclass}. By inspection of Tables \ref{eigIIexc} and \ref{eigIIclass} and Theorem \ref{cond3} we immediately conclude the following:

\blem \label{1testeII}The only bisymmetric triples satisfying
condition (\ref{condIIb}) are the cases \ref{cpdn7}, \ref{cpcn7},
\ref{cpg23}, \ref{cpe89}, \ref{cpe75} and

$$\bar{l}\xstrut\ref{cpan3}\, for\, $p=2l$,\, $n-p=2s$,\\

\xstrut\ref{cpdn4}\, for\, $p=2l$,\, $n-p=2s$,\\

\xstrut\ref{cpdn8}\, for\, $p=2l$,\\

\xstrut\ref{cpcn4}\, for\, $p=2l$,\, $n-p=2s$,\\

\xstrut\ref{cpcn8}\, for\, $p=2l$,\\

\xstrut\ref{cpe65}\, for\, $p=1$.\ear$$

For all other bisymmetric triples of Type II we can conclude that
there exists no Einstein adapted metric on $M$. \elem

Note that for all the triples listed in Lemma \ref{1testeII}, not only condition (\ref{condIIb}) is satisfied, but also $C_{\mfp_1}$ and $C_{\mfp_2}$ are scalar on $\mfn$. Thus $C_{\mfp}$ is scalar on $\mfn$ as well. We write

\beq \label{expcp1}C_{\mfp_a}\mid_{\mfn}=b_aId_{\mfn},\, a=1,2\textrm{ and }C_{\mfp}\mid_{\mfn}=bId_{\mfn},\textrm{ for
}b=b_1+b_2,\eeq

following the notation used in previous sections. We recall that $\ga_1$, $\ga_2$ are the eigenvalues of the Casimir operator $C_{\mfk}$ on $\mfp_1$, $\mfp_2$, respectively (see (\ref{gasdef})). Their values are also given in Tables \ref{eigIIexc} and \ref{eigIIclass}.
\li

\textbf{\emph{Proof of Theorem \ref{biII}:}} Einstein binormal  metrics are given by
Corollary \ref{binormal5}. First we observe that in order to exist
an Einstein binormal metric on $M$, $C_{\mfp}$ must be scalar on
$\mfn$ and $C_{\mfk}$ must be scalar on $\mfp$. The triples which
satisfy the first condition are those listed in Lemma
\ref{1testeII}. Furthermore, $C_{\mfk}$ is scalar on $\mfp$ if and only if $\ga_2=\ga_1$. We conclude
from Tables \ref{eigIIexc} and \ref{eigIIclass} that the spaces from Lemma \ref{1testeII}
which satisfy the condition $\ga_2=\ga_1$ are those listed below:

$$\bar{l|ccc} & \ga_1=\ga_2 & b_1 & b_2 \\ \hline
\xstrut\ref{cpan3},\textrm{ for }s=l=2p,\,n=4l &  \frac{1}{2} & \frac{1}{8} & \frac{1}{8}\\

\xstrut\ref{cpdn4}, \textrm{ for }s=l=2p,\,n=4l,\,l\geq 2 & \frac{2l-1}{4l-1} & \frac{l}{2(4l-1)} & \frac{l}{2(4l-1)}\\

\xstrut\ref{cpdn7}, \textrm{ for } n=2p,\, p\geq 2 & \frac{p-1}{2p-1} & \frac{p-1}{4(2p-1)} & \frac{p-1}{4(2p-1)}\\

\xstrut\ref{cpdn8}, \textrm{ for } p=2l,\, n=4l &  \frac{2l-1}{4l-1} & \frac{l}{2(4l-1)} & \frac{2l-1}{4(4l-1)}\\

\xstrut\ref{cpcn4}, \textrm{ for } s=l=2p,\,n=4l & \frac{2l+1}{4l+1} & \frac{l}{4(4l+1)} & \frac{l}{4(4l+1)}\\

\xstrut\ref{cpcn7}, \textrm{ for } n=2p & \frac{p+1}{2p+1} & \frac{p+1}{4(2p+1)} & \frac{p+1}{4(2p+1)}\\

\xstrut\ref{cpcn8}, \textrm{ for } p=2l,\,n=4l & \frac{2l+1}{4l+1} & \frac{l}{4(4l+1)} & \frac{2l+1}{4(4l+1)} \ear$$

In particular, we observe that for this list of spaces, $G$ is classical Lie group. For these spaces, Einstein binormal metrics, if exist, are given by positive solutions of
(\ref{polbin5}) in Corollary \ref{binormal5}. Since $c_{\mfk,\mfn}=\frac{1}{2}$, we simplify (\ref{polbin5}) and conclude that there exists an Einstein binormal metric if and only if

\beq\label{tII1}\triangle=1-2\ga(1-\ga+2b)\geq 0,\eeq

where $\ga=\ga_1=\ga_2$ and $b$ is the eigenvalue of $C_{\mfp}$ on
$\mfn$, i.e., $b=b_1+b_2$. In this case, the Einstein adapted metrics are
given by homotheties of

$$g_{\mfm}=B_{\mfp}\oplus XB_{\mfn},\textrm{
where }X=\frac{1\pm\sqrt{\triangle}}{2\ga}.$$

The values of $\triangle$ are as follows:

$$\bar{lc|lc} & \triangle & & \triangle\\ \hline

\xstrut\ref{cpan3} & 0 & \ref{cpcn4} &\frac{4l^2+2l+1}{(4l+1)^2}>0\\

\xstrut\ref{cpdn4} & \frac{1}{(4l-1)^2}>0 & \ref{cpcn7} & -\frac{1}{2p+1}<0\\

\xstrut\ref{cpdn7} & \frac{1}{2p-1}>0 & \ref{cpcn8} & \frac{l(2l-1)}{(4l+1)^2}>0\\

\xstrut\ref{cpdn8} & \frac{2l}{(4l-1)^2}>0 & &
\ear$$

Except in the case \ref{cpcn7}, there exists an Einstein adapted
metric. The values for $X$ are indicated in Table \ref{bimII}. Since the base space $N$ is isotropy irreducible, $g_N$ is Einstein. Moreover, it follows from Corollary \ref{binormal6} that $g_F$ is Einstein.

$\Box$

\li

\textbf{\emph{Proof of Theorem \ref{gfeinII}:} }The only cases which may admit an Einstein adapted metric $g_M$ are those listed in Lemma \ref{1testeII}, since they are the only
that satisfy the necessary condition (\ref{condIIb}). For bisymmetric fibrations of type II, Einstein adapted metrics $g_M$ such that $g_F$ is also Einstein are given by Corollary \ref{gfeincor2}. In particular, one of the two conditions in Corollary \ref{gfeincor} must be satisfied. In the first case, $\ga=\ga_1=\ga_2$, the Einstein adapted metrics with this property are binormal. Einstein binormal metrics are given by Theorem \ref{biII} and listed in Table
\ref{bimII}.
In the second case, we have $\ga_2=1-\ga_1$ and these metrics are described by (ii) in Corollary \ref{gfeincor2}. This case is possible only for the triple \ref{cpan3}
when $p=2l$ and $n-p=2s$. In this case $\ga_1=\frac{l}{l+s}$,
$\ga_2=\frac{s}{l+s}=1-\ga_1$, $b_1=\frac{\ga_1}{4}$ and
$b_2=\frac{\ga_2}{4}=\frac{1-\ga_1}{4}$. Note that on the formula for $D$, in Corollary \ref{gfeincor2}, $r=\frac{1}{2}$ since $N$ is a symmetric space. Using Corollary
\ref{gfeincor2}, we obtain $D(\ga_1)=0$ and thus $X_1=\frac{l+s}{2l}$, $X_2=\frac{l+s}{2s}$.

$\Box$

\li

\textbf{ \emph{Proof of Theorem \ref{genIIexc}:}} It follows from Lemma \ref{1testeII} that the only cases of Type II with $\mfg$ exceptional wich may admit an Einstein adapted metric are the cases \ref{cpe65} for $p=1$, \ref{cpg23}, \ref{cpe75} and \ref{cpe89}. We recall that for each of these spaces any Einstein adapted metric is of the form

$$g_{\mfm}=\frac{1}{X_1}B_{\mfp_1}\oplus\frac{1}{X_2}B_{\mfp_2}\oplus B_{\mfn},$$ where
$X_1$ and $X_2$ are positive solutions of the system of equations given in Theorem \ref{eqFsym} which are as follows:

\beqar
\label{eqII112}2\ga_1X_1^2X_2+(1-\ga_1)X_2-2\ga_2X_1X_2^2-(1-\ga_2)X_1=0,
\\
\label{eqII122}2b_1X_2+2b_2X_1-2X_1X_2+2\ga_1X_1^2X_2+(1-\ga_1)X_2=0.
\eeqar

Also we recall that the eigenvalues $b_i$ and $\ga_i$, for $i=1,2$, can be found in Table \ref{eigIIexc}. We explain briefly how to obtain solutions of these equations for each case. Most of the calculations are performed with Maple and we omit the details. This proof can be found in more detail in \cite[\S 3.4]{Fa1}. For each case, the solution is given as follows: $X_1=\al$, $X_2=f(\al)$, where $f$ is a polynomial with degree 3 and $\al$ is a root of a polynomial $t$ with degree $4$. Hence, the number of Einstein adapted metrics is the number of positive roots $\al$ of $t$ such that $f(\al)>0$. We indicate the polynomials $f$ and $t$ for each case and the number of Einstein adapted metrics. Approximations of the solutions are given in Table \ref{tabgenII}. For exact solutions see \cite[\S 3.4]{Fa1}.

$$\bar{ll}\ref{cpg23}) & f(\al)=-\frac{7}{4}\al^3+12\al^2-\frac{5899}{36}\al+19\\
 & t(z) = 63z^4-432z^3+1088z^2-1224z+513\\
 & 2 \textrm{ Einstein adapted metrics} \\ \\

 \ref{cpe65}) & f(\al)=-\frac{156}{7}\al^3+\frac{552}{7}\al^2-\frac{571}{7}\al+\frac{176}{7}\\
& t(z)=234z^4-828z^3+993z^2-474z+77\\
& 4 \textrm{ Einstein adapted metrics} \\ \\

\ref{cpe75}, \textrm{ for }$p=2$)  & f(\al)=-\frac{140}{3}\al^3+148\al^2-\frac{681}{5}\al+\frac{184}{5}\\
& t(z)=350z^4-1110z^3+1179z^2-492z+69\\
& 4 \textrm{ Einstein adapted metrics} \\ \\

\ref{cpe75}, \textrm{ for }$p=4$) & f(\al)=-\frac{100}{3}\al_i^3+100\al_i^2-\frac{262}{3}\al_i+26\\
& t(z)=200z^4-600z^3+614z^2-264z+39\\
 & 2 \textrm{ Einstein adapted metrics} \\ \\

 \ref{cpe75}, \textrm{ for } $p=6$) &  f(\al)=-\frac{2500}{3}z^3+820z^2-235z+24\\
& t(z)=1250z^4-1230z^3+415z^2-60z+3\\
 & 2 \textrm{ Einstein adapted metrics}\\ \\

\ref{cpe89}) & t(z)=9z^4-195z^3+1198z^2-1395z+464\textrm{ has no positive roots }\\
 & \textrm{No Einstein adapted metrics} \ear$$

$\Box$

\li

\textbf{\emph{Proof of Theorem \ref{genIIall}:}} This an application of Corollaries \ref{contcor1} and \ref{contcor2}. The only case when $\ga_2=1-\ga_1$ is \ref{cpan3}, for
$p=2l$ and $n=2(l+s)$. For this case, we have $\ga_1=\frac{l}{l+s}$,
$\ga_2=\frac{s}{l+s}$, $b_1=\frac{l}{4(l+s)}=\frac{\ga_1}{4}$ and
$b_2=\frac{s}{4(l+s)}=\frac{1-\ga_1}{4}$. Since the restriction to the fiber $g_F$ is not Einstein, the required metric
is given by Corollary \ref{contcor2} (ii). A simple
calculation show that $D(\ga_1)=-\frac{1}{2}\ga_1(1-\ga_1)<0$, for
every $l,s$, since $0<\ga_1<1$. Hence, in this case there are no
other Einstein adapted metrics besides those found previously.

The cases such that $\ga_2=\ga_1$ are those listed in the proof of
Theorem \ref{biII}. In this case, there exists an Einstein adapted
metric such that $g_F$ is not Einstein if and only if
$D(\ga_1)\geq 0$, where

$$D(\ga_1)=4r^2(1-\ga_1)-2\ga_1(2b_2+1-\ga_1)(2b_1+1-\ga_1),$$

according to Corollary \ref{contcor1} (ii). Since $N=G/K$ is a symmetric space, $r=\frac{1}{2}$.

\li

$$\bar{ll}\ref{cpan3},\,s=l=2p,\,n=4l) & \ga_1=\frac{1}{2}, \, b_1=b_2=\frac{1}{8} \\
& D(\ga_1)=\frac{1}{2}(-\ga_1^3+4\ga_1^2-6\ga_1+2)<0\\
& \textrm{No metric}\\ \\

\ref{cpdn7},\, n=2p,\, p\geq 2 & \ga_1=\frac{p-1}{2p-1} , \, b_1=b_2=\frac{\ga_1}{4}\\
& D(\ga_1)=\frac{1}{2}(-\ga_1^3+4\ga_1^2-6\ga_1+2)>0\textrm{ iff }p=2,\ldots,6\\
& 2\textrm{ metrics for }p=2,\ldots,6\\ \\

\ref{cpcn7},\,n=2p) & \ga_1=\frac{p+1}{2p+1} , \, b_1=b_2=\frac{\ga_1}{4}\\
& D(\ga_1)=\frac{1}{2}(-\ga_1^3+4\ga_1^2-6\ga_1+2)<0, \forall p\\
& \textrm{No metric}\\ \\

\ref{cpdn4},\,s=l=2p\geq 2,\, n=4l) & \ga_1=\frac{2l-1}{4l-1},\,b_1=b_2=\frac{l}{2(4l-1)}=\frac{1-\ga_1}{4}\\
& D(\ga_1)=-\frac{1}{2}(\ga_1-1)(3\ga_1-1)(3\ga_1-2)<0,\,\forall l\\
& \textrm{No metric}\\ \\

\ref{cpdn8},\,p=2l,\,n=4l) & \ga_1=\frac{2l-1}{4l-1},\,b_1=\frac{l}{2(4l-1)}=\frac{1-\ga_1}{4},\,b_2=\frac{2l-1}{4(4l-1)}=\frac{\ga_1}{4}\\
& D(\ga_1)=\frac{1}{2}(1-\ga_1)(3\ga_1^2-6\ga_1+2)>0\textrm{ iff }l=1\\
& 2\textrm{ metrics for }l=1\\ \\

\ref{cpcn4},\,s=l=2p,\,n=4l) & \ga_1=\frac{2l+1}{4l+1},\,b_1=b_2=\frac{l}{4(4l+1)}=\frac{1-\ga_1}{8}\\
& D(\ga_1)=\frac{1}{8}(1-\ga_1)(25\ga_1^2-25\ga_1+8)>0, \forall l\\
& 2\textrm{ metrics for }, \forall l\\ \\

\ref{cpcn8},\,p=2l,\,n=4l) & \ga_1=\frac{2l+1}{4l+1},\,b_1=\frac{l}{4(4l+1)}=\frac{1-\ga_1}{8},\,b_2=\frac{2l+1}{4(4l+1)}=\frac{\ga_1}{4}\\
& D(\ga_1)=\frac{1}{4}(1-\ga_1)(5\ga_1^2-10\ga_1+4)>0, \textrm{ iff }l\geq 3\\
& 2\textrm{ metrics for }l\geq 3
\ear$$

\li

For each case, we calculate $X_1$, $X_2$ according to the formulas given by Corollaries \ref{contcor1} and \ref{contcor2}. These values are indicated in Table \ref{nonbimII}.

$\Box$

\section{Application to $4$-symmetric Spaces}\label{4symsection}

A homogeneous space $G/L$ is said to be a $4$-symmetric space if
there exists $\sg\in Aut(G)$ such that

$$(G_{\sg})_0\subset L\subset G_{\sg}$$

and $\sg$ has order $4$. Compact simply connected irreducible
$4$-symmetric spaces have been classified by J.A.Jimenez in
\cite{Ji} following the previous work of V.Ka\v{c} (see e.g.
\cite{He}, Chap.X), J.A Wolf and A.Gray \cite{WG}. It is shown in
\cite{Ji} that any compact simply connected irreducible $4$-symmetric space is the
total space of a fiber bundle whose fiber and base space are
symmetric spaces and the base is an isotropy irreducible space of
maximal rank. These spaces are fully described in Tables III, IV and
V in \cite{Ji}. Hence, for each compact simply connected irreducible
$4$-symmetric space $M=G/L$ there is a bisymmetric fibration
$M=G/L \rightarrow G/K = N$ such that $N$ is isotropy irreducible and $K$ has maximal rank. Therefore, we can apply the results obtained in previous sections. We recall that the bisymmetric fibrations $M=G/L \rightarrow G/K = N$ considered in this paper are such that $L$ has maximal rank. This allow us to take conclusions for $4$-symmetric spaces of maximal rank. The bisymmetric triples $(\mfg,\mfk,\mfl)$
corresponding to $4$-symmetric spaces of maximal rank must be some
of Tables \ref{eigIexc}, \ref{eigIclass}, \ref{eigIIexc} and
\ref{eigIIclass}. Each of these tables contains a column which indicates if the space is a $4$-symmetric space. Hence, a simple comparison between the Tables in this paper and the classification in \cite{Ji} allow us to easily
conclude about the existence of Einstein metrics on $4$-symmetric
spaces.

\bcor\label{4sym1} Let $M=G/L\rightarrow G/K$ be a bisymmetric fibration such that $M$ is a compact simply-connected irreducible $4$-symmetric space (of maximal rank). If $M$ admits an Einstein adapted metric, then the corresponding bisymmetric triple $(\mfg,\mfk,\mfl)$ is one of the triples in Tables \ref{mIexc}, \ref{mIclass}, \ref{bimII}, \ref{nonbimII} and \ref{tabgenII}.
\ecor

\newpage

\section{Tables}\label{alltabs}

\btab[h]\caption{Dual Coxeter Numbers}\label{tabcoxeter}

$$\bar{|c|c|}

\hline \textrm{Coxeter group} & \textrm{Dual Coxeter number}\\\hline

A_n & n+1\\

B_n & 2n-1\\

C_n & n+1\\

D_n & 2n-2\\

E_6 & 12\\

E_7 &  18\\

E_8 & 30\\

F_4 & 9\\

G_2 & 4\\ \hline

\ear$$ \etab

\btab[h]\caption{Symmetric pairs of compact type of maximal rank -
Exceptional Spaces}\label{spexc}
$$\bar{|l|l||l|l|} \hline \mfg & \mfk & \mfg & \mfk
\\\hline

\xstrut\mff_4 & \mfsp_3\oplus\mfsu_2 & \mfe_7 & \mfsu_8\\

\xstrut\mff_4 & \mfso_9 & \mfe_7 & \mfe_6\oplus\reals\\

\xstrut\mfg_2 & \mfsu_2\oplus \mfsu_2 & \mfe_7 & \mfso_{12}\oplus \mfsu_2\\

\xstrut\mfe_6 & \mfso_{10}\oplus \reals & \mfe_8 & \mfso_{16}\\

\xstrut\mfe_6 & \mfsu_6\oplus\mfsu_2 & \mfe_8 & \mfe_7\oplus
\mfsu_2\\ \hline\ear$$\etab

\btab[h]\caption{Symmetric pairs of compact type of maximal rank -
Classical Spaces}\label{spclass}
$$\bar{|l|l|} \hline \mfg & \mfk
\\\hline

\xstrut\mfso_{2n} & \mfu_n\\

\xstrut\mfso_n & \mfso_{2p}\oplus\mfso_{n-2p}\\

\xstrut\mfsp_{n} & \mfu_n\\

\xstrut\mfsp_{n} & \mfsp_p\oplus\mfsp_{n-p}\\

\xstrut\mfsu_n &
\mfsu_p\oplus\mfsu_{n-p}\oplus\reals\\\hline\ear$$\etab

\btab[c]\caption{Bisymmetric triples of type I and their eigenvalues
- Exceptional spaces}\label{eigIexc}
$$\bar{|c|lll|c|c|}\hline \xstrut  & \mfg & \mfk  & \mfl & \ga & b^j\\ \hline

\xstrut\ref{cpf41} & \mff_4 & \mfso_9 &
\mfso_p\oplus\mfso_{9-p},\,p=1,3,5,7 & \frac{7}{9} &
\frac{p(9-p)}{72}\\ \hline

\xstrut\ref{cpf42} &\mff_4 & \mfsp_3\oplus\mfsu_2 &
\mfsp_3\oplus\reals &
\frac{2}{9} & \frac{1}{18}\\

\xstrut\ref{cpf43} & & & \mfu_3\oplus\mfsu_2 & \frac{4}{9} &
\frac{1}{4},\,\frac{2}{9}\\

\xstrut\ref{cpf44} & & & \mfsp_2\oplus\mfsu_2\oplus \mfsu_2 &
\frac{4}{9} & \frac{1}{9},\,\frac{1}{18}\\\hline

\xstrut\ref{cpg21} & \mfg_2 & \mfsu_2\oplus \mfsu_2 &
\reals\oplus\mfsu_2 &
\frac{1}{2} & \frac{1}{8}\\

\xstrut\ref{cpg22} & &  & \mfsu_2\oplus\reals  & \frac{1}{6} &
\frac{1}{6}\\\hline

\xstrut\ref{cpe81} & \mfe_8 & \mfso_{16} &
\mfso_{2p}\oplus\mfso_{16-2p},\,p=1,\ldots,4 & \frac{1}{5} &
\frac{p(8-p)}{60}\\

\xstrut\ref{cpe82} & & & \mfu_8 & \frac{1}{5} &
\frac{4}{15},\,\frac{3}{15},\,\frac{7}{15} \\\hline

\xstrut\ref{cpe83} & \mfe_8 & \mfe_7\oplus \mfsu_2 & \mfe_7\oplus \reals & \frac{1}{15} & \frac{1}{60}\\

\xstrut\ref{cpe84} & & & \mfe_6\oplus\reals\oplus \mfsu_2 &
\frac{3}{5} &
\frac{11}{60},\frac{9}{20}\\

\xstrut\ref{cpe86} & & & \mfso_{12}\oplus\mfsu_2\oplus\mfsu_2 &
\frac{3}{5}
& \frac{4}{15},\frac{1}{5}\\

\xstrut\ref{cpe88} & & & \mfsu_8\oplus \mfsu_2 & \frac{3}{5} &
\frac{1}{4}\\\hline

\xstrut\ref{cpe71} & \mfe_7 & \mfso_{12}\oplus \mfsu_2 &
\mfso_{12}\oplus\reals & \frac{1}{9} & \frac{1}{36}\\

\xstrut\ref{cpe72} & & & \mfu_6\oplus \mfsu_2 & \frac{5}{9} &
\frac{1}{6},\,\frac{5}{18}\\

\xstrut\ref{cpe74} & & & \mfso_p\oplus\mfso_{12-p}\oplus \mfsu_2,
\,p=2,4,6 & \frac{5}{9} &  \frac{p(12-p)}{144}\\\hline

\xstrut\ref{cpe76} & \mfe_7 & \mfe_6\oplus\reals &
\mfso_{10}\oplus\reals\oplus\reals & \frac{2}{3}
&\frac{2}{9},\frac{1}{6},\,\frac{4}{9}\\

\xstrut\ref{cpe77} &  &  & \mfsu_6\oplus\mfsu_2\oplus\reals &
\frac{2}{3} & \frac{5}{18},\frac{2}{9}\\\hline

\xstrut\ref{cpe78} & \mfe_7 & \mfsu_8 &
\mfsu_p\oplus\mfsu_{8-p}\oplus\reals,\,1\leq p\leq 4 & \frac{4}{9} &
\bar{cc} \xstrut\frac{1}{9}, & p=1\\ \xstrut\frac{2}{9},\,\frac{1}{6}, & p=2\\
\xstrut\frac{2}{9},\,\frac{1}{3}, & p=3\\ \xstrut\frac{2}{9},
\frac{4}{9},\,\frac{11}{36}, & p=4 \ear
\\\hline

\xstrut\ref{cpe61} & \mfe_6 & \mfso_{10}\oplus \reals &
\mfu_5\oplus\reals
& \frac{2}{3} & \frac{5}{12},\,\frac{1}{6},\,\frac{1}{4}\\

\xstrut\ref{cpe62} & & & \mfso_p\oplus\mfso_{10-p}\oplus\reals,
\,p=2,4 & \frac{2}{3} & \frac{p(10-p)}{96}
\\\hline

\xstrut\ref{cpe63} & \mfe_6 & \mfsu_6\oplus\mfsu_2 &
\mfsu_6\oplus\reals &
\frac{1}{6} & \frac{1}{24}\\

\xstrut\ref{cpe64} & & &
\mfsu_p\oplus\mfsu_{6-p}\oplus\reals\oplus\mfsu_2 & \frac{1}{2} &
\frac{p+2}{24},\,\frac{p}{8}\\ \hline \ear$$\etab

\li

\bland\btab[c]\caption{Bisymmetric triples of type I and their
eigenvalues - Classical spaces}\label{eigIclass}
$$\bar{|c|lll|c|c|}\hline  & \mfg & \mfk  & \mfl & \ga & b^j\\ \hline

\xstrut\ref{cpan1} & \mfsu_n & \mfsu_p\oplus\mfsu_{n-p}\oplus\reals
& \mfsu_l\oplus\mfsu_{p-l}\oplus\reals\oplus\mfsu_{n-p}\oplus\reals
& \frac{p}{n} & \frac{p-l}{2n},\,\frac{l}{2n}\\ \hline

\xstrut\ref{cpbn1} & \mfso_{2n+1} &
\mfso_{2p+1}\oplus\mfso_{2(n-p)},\,p=0,\ldots,n-1 &
\mfso_{2l+1}\oplus\mfso_{2(p-l)}\oplus\mfso_{2(n-p)} &
\frac{2p-1}{2n-1} & \frac{p-l}{2n-1},\,\frac{4l+1}{4(2n-1)}\\

\xstrut\ref{cpbn2} & & &
\mfso_{2p+1}\oplus\mfso_{2s}\oplus\mfso_{2(n-p-s)}
& \frac{2(n-p-1)}{2n-1} & \frac{n-p-s}{2n-1},\,\frac{s}{2n-1}\\

\xstrut\ref{cpbn3} & & & \mfso_{2p+1}\oplus\mfu_{n-p} &
\frac{2(n-p-1)}{2n-1} & \frac{n-p-1}{2(2n-1)}\\ \hline

\xstrut\ref{cpdn1} & \mfso_{2n} & \mfu_n &
\mfu_p\oplus\mfu_{n-p},\,p=1,\ldots,n-1 & \frac{n}{2(n-1)} &
\frac{n-p}{2(n-1)},\,\frac{p}{2(n-1)},\,\frac{n-2}{4(n-1)}\\ \hline

\xstrut\ref{cpdn2} & \mfso_{2n} &
\mfso_{2p}\oplus\mfso_{2(n-p)},\,p=1,\ldots,\lfloor
\frac{n}{2}\rfloor  &
\mfso_{2l}\oplus\mfso_{2(p-l)}\oplus\mfso_{2(n-p)}& \frac{p-1}{n-1}
& \frac{p-l}{2(n-1)},\,\frac{l}{2(n-1)}\\

\xstrut\ref{cpdn5} & & & \mfu_p\oplus\mfso_{2(n-p)} &
\frac{p-1}{n-1} & \frac{p-1}{4(n-1)}\\ \hline

\xstrut\ref{cpcn1} & \mfsp_{n} & \mfu_n &
\mfu_p\oplus\mfu_{n-p},\,p=1,\ldots,n-1 & \frac{n}{2(n+1)} &
\frac{n-p}{2(n+1)},\,\frac{n-p}{n+1},\,\frac{p}{2(n+1)},\,\frac{p}{n+1},\,\frac{n+2}{2(n+1)}
 \\ \hline

\xstrut\ref{cpcn2} & \mfsp_n & \mfsp_p\oplus\mfsp_{n-p} &
\mfsp_l\oplus\mfsp_{p-l}\oplus\mfsp_{n-p} & \frac{p+1}{n+1} &
\frac{p-l}{4(n+1)},\,\frac{l}{4(n+1)}\\

\xstrut\ref{cpcn5} &  &  & \mfu_p\oplus\mfsp_{n-p} & \frac{p+1}{n+1}
& \frac{p+1}{4(n+1)}\\ \hline \ear$$\etab \eland

\bland\btab[c]\caption{Bisymmetric triples of type II and their
eigenvalues - Exceptional
spaces}\label{eigIIexc}$$\bar{|c|lll|c|c|c|c|}\hline  &\mfg &
\mfk  & \mfl & \ga_1 & \ga_2 & b_1^{\phi} & b_2^{\phi} \\\hline

\xstrut\ref{cpf45} & \mff_4 & \mfsp_3\oplus\mfsu_2 & \mfu_3\oplus
\reals & \frac{4}{9} & \frac{2}{9} &
\big(\frac{1}{4},\frac{2}{9}\big) &
\frac{1}{18}\\

\xstrut\ref{cpf46} & & & \mfsu_2\oplus\mfsp_2\oplus\reals &
\frac{4}{9} &
\frac{2}{9} & \big(\frac{1}{9},\frac{1}{18}\big) & \frac{1}{18}\ \\
\hline

\xstrut\ref{cpg23} & \mfg_2 & \mfsu_2\oplus \mfsu_2  &
\reals\oplus\reals & \frac{1}{2} & \frac{1}{6} & \frac{1}{8} &
\frac{1}{6}\\\hline

\xstrut\ref{cpe85} & \mfe_8 & \mfe_7\oplus \mfsu_2 &
\mfe_6\oplus\reals\oplus\reals & \frac{3}{5} & \frac{1}{15} &
\big(\frac{11}{60},\frac{11}{60},\frac{9}{20}\big) & \frac{1}{60}\\

\xstrut\ref{cpe87} & & & \mfso_{12}\oplus\mfsu_2\oplus\reals &
\frac{3}{5} & \frac{1}{15} & \big(\frac{4}{15},\frac{1}{5}\big)  &
\frac{1}{60}
\\

\xstrut\ref{cpe89} & & & \mfsu_8\oplus\reals & \frac{3}{5} &
\frac{1}{15} & \frac{1}{4}& \frac{1}{60} \\ \hline

\xstrut\ref{cpe73} & \mfe_7 & \mfso_{12}\oplus \mfsu_2 &
\mfu_6\oplus \reals & \frac{5}{9} & \frac{1}{9} &
\big(\frac{5}{18},\frac{1}{6},\frac{5}{18}\big) & \frac{1}{36}\\

\xstrut\ref{cpe75} & & & \mfso_p\oplus\mfso_{12-p}\oplus \reals,
\,p\textrm { even} & \frac{5}{9} & \frac{1}{9} & \frac{1}{36}&
\frac{p(12-p)}{144}\\ \hline

\xstrut\ref{cpe65} & \mfe_6 & \mfsu_6\oplus\mfsu_2 &
\mfsu_p\oplus\mfsu_{6-p}\oplus \reals\oplus\reals & \frac{1}{2} &
\frac{1}{6}
 & \frac{1}{24} & \frac{p+2}{24},\, \frac{p}{8}\\ \hline
\ear$$\etab\eland

\bland \btab[c]\caption{Bisymmetric triples of type II and their
eigenvalues -
Classical spaces}\label{eigIIclass}$$\bar{|c|lll|c|c|c|c|}\hline  & \mfg & \mfk  & \mfl & \ga_1 & \ga_2 & b_1^{\phi} & b_2^{\phi} \\
\hline
\xstrut\ref{cpan3} & \mfsu_n & \mfsu_p\oplus\mfsu_{n-p}\oplus\reals
&
\mfsu_l\oplus\mfsu_{p-l}\oplus\mfsu_s\oplus\mfsu_{n-p-s}\oplus\reals\oplus\reals\oplus\reals
& \frac{p}{n} & \frac{n-p}{n} &
\big(\frac{p-l}{2n},\frac{l}{2n}\big) &
\big(\frac{n-p-s}{2n},\frac{s}{2n}\big) \\ \hline

\xstrut\ref{cpbn5} & \mfso_{2n+1} & \mfso_{2p+1}\oplus
\mfso_{2(n-p)} & \mfso_{2l+1}\oplus
\mfso_{2(p-l)}\oplus\mfso_{2s}\oplus\mfso_{2(n-p-s)} &
\frac{2p-1}{2n-1} & \frac{2(n-p-1)}{2n-1} &
\big(\frac{p-l}{2n-1},\frac{4l+1}{4(2n-1)}\big) &
\big(\frac{n-p-s}{2n-1},\frac{s}{2n-1}\big)\\

\xstrut\ref{cpbn4} &  & & \mfso_{2l+1}\oplus
\mfso_{2(p-l)}\oplus\mfu_{n-p} & \frac{2p-1}{2n-1} &
\frac{2(n-p-1)}{2n-1}  &
\big(\frac{p-l}{2n-1},\frac{4l+1}{4(2n-1)}\big) &
\frac{n-p-1}{2(2n-1)}\\ \hline

\xstrut\ref{cpdn4} & \mfso_{2n} & \mfso_{2p}\oplus \mfso_{2(n-p)} &
\mfso_{2l}\oplus
\mfso_{2(p-l)}\oplus\mfso_{2s}\oplus\mfso_{2(n-p-s)} &
\frac{p-1}{n-1}& \frac{n-p-1}{n-1} &
\big(\frac{p-l}{2(n-1)},\frac{l}{2(n-1)}\big) &
\big(\frac{n-p-s}{2(n-1)},\frac{s}{2(n-1)}\big)\\

\xstrut\ref{cpdn7} &  & & \mfu_p\oplus \mfu_{n-p} & \frac{p-1}{n-1}&
\frac{n-p-1}{n-1} & \frac{p-1}{4(n-1)} & \frac{n-p-1}{4(n-1)}\\

\xstrut\ref{cpdn8} &  & & \mfso_{2l}\oplus
\mfso_{2(p-l)}\oplus\mfu_{n-p} & \frac{p-1}{n-1}& \frac{n-p-1}{n-1}
& \big(\frac{p-l}{2(n-1)},\frac{l}{2(n-1)}\big) &
\frac{n-p-1}{4(n-1)}\\ \hline

\xstrut\ref{cpcn4} & \mfsp_n & \mfsp_p\oplus \mfsp_{n-p} &
\mfsp_l\oplus \mfsp_{p-l}\oplus\mfsp_{s}\oplus\mfsp_{n-p-s} &
\frac{p+1}{n+1}& \frac{n-p+1}{n+1} &
\big(\frac{p-l}{4(n+1)},\frac{l}{4(n+1)}\big) &
 \big(\frac{n-p-s}{4(n+1)},\frac{s}{4(n+1)}\big)\\

\xstrut\ref{cpcn7} &  &  & \mfu_p\oplus \mfu_{n-p} &
\frac{p+1}{n+1}&
\frac{n-p+1}{n+1} & \frac{p+1}{4(n+1)}& \frac{n-p+1}{4(n+1)}\\

\xstrut\ref{cpcn8} &  &  & \mfsp_l\oplus \mfsp_{p-l}\oplus\mfu_{n-p}
& \frac{p+1}{n+1}& \frac{n-p+1}{n+1} &
\big(\frac{p-l}{4(n+1)},\frac{l}{4(n+1)}\big) &
\frac{n-p+1}{4(n+1)}\\ \hline

 \ear$$\etab \eland

\btab[c]
\begin{minipage}{\textwidth}\caption{Einstein Bisymmetric triples of type I -
Exceptional spaces}\label{mIexc}
$$\bar{|lll|c|c|}\hline\mfg & \mfk  & \mfl & 4\textrm{-sym.} & X\\ \hline

\xstrut\mff_4 & \mfsp_3\oplus\mfsu_2 & \mfsp_3\oplus\reals & \textrm{ yes} &
\frac{1}{2},\,4\,^{(a)}\\\hline

\xstrut\mff_4 & \mfso_9 & \mfso_8 & \textrm{no} & 1,\,\frac{2}{7}\,^{(c)}\\

\xstrut& & \mfso_7\oplus\reals & \textrm{ yes}& \frac{9\pm\sqrt{8}}{14}\,^{(a)}\\\hline

\xstrut\mfg_2 & \mfsu_2\oplus \mfsu_2 & \reals\oplus\mfsu_2 & \textrm{ yes} &
\frac{1}{2},\,\frac{3}{2}\\

\xstrut&  & \mfsu_2\oplus\reals & \textrm{ yes} & \frac{6\pm\sqrt{22}}{2}\\\hline

\xstrut\mfe_6 & \mfso_{10}\oplus \reals &
\mfso_{8}\oplus\reals\oplus\reals & \textrm{ yes} &  1,\,\frac{1}{2}\,^{(c)}\\\hline

\xstrut\mfe_6 & \mfsu_6\oplus\mfsu_2 & \reals\oplus\mfsu_5\oplus\mfsu_2 & \textrm{ yes} & \frac{1}{2},\,\frac{3}{2}\,^{(a)}\\

\xstrut & & \mfsu_6\oplus\reals & \textrm{ yes} &
\frac{1}{2},\,\frac{11}{2}\,^{(a)}\\\hline

\xstrut\mfe_7 & \mfso_{12}\oplus \mfsu_2 & \mfso_{12}\oplus\reals & \textrm{ yes} & \frac{17}{2},\,\frac{1}{2}\\

\xstrut & & \reals\oplus\mfso_{10}\oplus \mfsu_2 & \textrm{no} &
\frac{1}{2},\,\frac{13}{10}\\

\xstrut & & \mfso_4\oplus\mfso_{8}\oplus \mfsu_2  & \textrm{no}  & 1,\, \frac{4}{5} \,^{(c)}\\\hline

\xstrut\mfe_7 & \mfsu_8 & \mfsu_7\oplus\reals & \textrm{no} &
\frac{1}{2},\,\frac{7}{4}\\\hline

\xstrut\mfe_8 & \mfe_7\oplus \mfsu_2 & \mfe_7\oplus \reals & \textrm{no} &
\frac{1}{2},\,\frac{29}{2}\,^{(a)}\\\hline

\xstrut\mfe_8 & \mfso_{16} & \mfso_{2p}\oplus\mfso_{16-2p} & \textrm{yes} &
 \frac{15\pm\sqrt{7p^2-56p+113}}{14} \,^{(b)} \\ \hline
 \ear$$

\footnotetext{See Theorem \ref{mtypeI}. Compare with Table \ref{eigIexc}.}
\footnotetext{$^{(a)}$ Metrics also obtained by Dickinson and Kerr in \cite{DK}.}
\footnotetext{$^{(b)}$ \emph{idem} for $p=1,3$. For $p=4$, one of the metrics is the standard metric obtained by Wang and Ziller in \cite{WZ}.}
\footnotetext{$^{(c)}$ The standard metric was obtained by Wang and Ziller in \cite{WZ}.}
\end{minipage}\etab

\bland\btab[c]\caption{Einstein Bisymmetric triples of type I -
Classical spaces}\label{mIclass}
\begin{minipage}{\textwidth}
$$\bar{|lll|c|c|}\hline\mfg & \mfk  & \mfl & 4\textrm{-sym.} & X\\ \hline

\xstrut\mfso_{2n} & \mfso_{2p}\oplus\mfso_{2(n-p)} &
\mfso_{p}\oplus\mfso_{p}\oplus\mfso_{2(n-p)},\,p\,even & \textrm{no} &
\frac{n-1\pm\sqrt{p^2-(2n+1)p+n^2+1}}{2(p-1)}\\

\xstrut& & \mfu_p\oplus\mfso_{2(n-p)} & \textrm{yes} &
\frac{1}{2},\,\frac{n}{p-1}-\frac{1}{2}\,^{(a)}\\\hline

\xstrut\mfso_{2n+1} & \mfso_{2p+1}\oplus\mfso_{2(n-p)} &
\mfso_{2p+1}\oplus\mfso_{n-p}\oplus\mfso_{n-p},\,n-p\,even & \textrm{yes} &
\frac{2n-1\pm\sqrt{4p^2+8p-4n+5}}{4(n-p-1)}\\

\xstrut& & \mfso_{2p+1}\oplus\mfu_{n-p} & \textrm{no} &
\frac{1}{2},\,\frac{n+p}{2(n-p-1)}\,^{(a)}\\\hline

\xstrut\mfsp_n & \mfsp_{2l}\oplus\mfsp_{n-2l} &
\mfsp_l\oplus\mfsp_l\oplus\mfsp_{n-2l} & \textrm{no} &
\frac{n+1\pm\sqrt{6l^2+(3-4n)l+n^2+1}}{2(2l+1)}\\

\xstrut & \mfsp_{p}\oplus\mfsp_{n-p} & \mfu_p\oplus\mfsp_{n-p} & \textrm{yes} & \frac{1}{2},\,\frac{1}{2}+\frac{n-p}{p+1}\,^{(a)}\\\hline

\xstrut\mfsu_n & \mfsu_{2l}\oplus\mfsu_{n-2l}\oplus\reals &
\mfsu_l\oplus\mfsu_l\oplus\reals\oplus\mfsu_{n-2l}\oplus\reals & \textrm{no}
 & \frac{1}{2},\,\frac{n}{2l}-\frac{1}{2}\\ \hline\ear$$

\footnotetext{See Theorem \ref{mtypeI}. Compare with Table \ref{eigIclass}.}
\footnotetext{$^{(a)}$ Metrics also obtained by Dickinson and Kerr in \cite{DK}.}
\end{minipage}
\etab\eland

\btab[h!]\caption{Bisymmetric triples of type II with Einstein
metric such that $g_F$ is also Einstein}\label{bimII}
\begin{minipage}{\textwidth}

\bc{$g_M$ binormal}\ec

$$\bar{|lll|c|c|}\hline\mfg & \mfk  & \mfl  & 4\textrm{-sym.} & X\\ \hline

\xstrut\mfsu_{4l} & \mfsu_{2l}\oplus\mfsu_{2l}\oplus\reals &
\mfsu_l\oplus\mfsu_l\oplus\oplus\mfsu_l\oplus\mfsu_l\oplus\reals^3 & \textrm{no} & 1\,\,^{(a)}\\

\xstrut\mfso_{8l} & \mfso_{4l}\oplus\mfso_{4l} & \mfso_{2l}
\oplus\mfso_{2l} \oplus\mfso_{2l} \oplus\mfso_{2l} & \textrm{yes} &
1,\,\frac{2l}{2l-1}\,\,^{(a)}\\

\xstrut\mfso_{8l} & \mfso_{4l}\oplus\mfso_{4l} & \mfso_{2l}
\oplus\mfso_{2l} \oplus\mfu_{2l} & \textrm{no} &
\frac{4l-1\pm\sqrt{2l}}{2(2l-1)}\\

\xstrut\mfso_{4l} & \mfso_{2l}\oplus\mfso_{2l} & \mfu_{l}
\oplus\mfu_{l},\,l\geq 2 & \textrm{no} & \frac{2l-1\pm\sqrt{2l-1}}{2(l-1)}\\

\xstrut\mfsp_{4l} & \mfsp_{2l}\oplus\mfsp_{2l} &
\mfsp_l\oplus\mfsp_l\oplus\mfsp_l\oplus\mfsp_l & \textrm{no} &
\frac{4l+1\pm\sqrt{4l^2+2l+1}}{2(2l+1)}\\

\xstrut\mfsp_{4l} & \mfsp_{2l}\oplus\mfsp_{2l} &
\mfsp_l\oplus\mfsp_l\oplus\mfu_{2l} & \textrm{no} &
\frac{4l+1\pm\sqrt{l(2l-1)}}{2(2l+1)}
\\ \hline\ear$$

\li

\bc{$g_M$ non-binormal}\ec

$$\bar{|lll|c|c|c|}\hline\mfg & \mfk  & \mfl & 4\textrm{-sym.} & X_1 & X_2\\ \hline

\xstrut\mfsu_{2(l+s)} & \mfsu_{2l}\oplus\mfsu_{2s}\oplus\reals &
\mfsu_l\oplus\mfsu_l\oplus\mfsu_s\oplus\mfsu_s\oplus\reals^3 & \textrm{yes} &
\frac{l+s}{2l} & \frac{l+s}{2s}\\\hline\ear$$

\footnotetext{See Theorems \ref{biII} and \ref{gfeinII}. Compare with Tables \ref{eigIIexc} and \ref{eigIIclass}.}
\footnotetext{$^{(a)}$ The standard metric was obtained by Wang and Ziller in \cite{WZ}.}
\end{minipage}\etab

\li

\li

\li

\btab[h!]\caption{All other Einstein adapted metrics for the
bisymmetric triples of Type II which admit an EAM $g_M$ such that
$g_F$ is Einstein}\label{nonbimII}
\begin{minipage}{\textwidth}
$$\bar{|lll|c|c|c|}\hline\mfg & \mfk  & \mfl & 4\textrm{-sym.} & X_1 & X_2\\ \hline

\xstrut\mfso_{4l}& \mfso_{2l}\oplus\mfso_{2l} &
\mfu_l\oplus\mfu_l,
\,l=2,\ldots,6 & \textrm{no}  & \frac{2l(l-1)\pm\sqrt{(-l^4+7l^3-5l^2+l)/2}}{2(l-1)(3l-1)} & \frac{l}{2(l-1)}.\frac{1}{X_1}\\

\xstrut\mfso_{8} & \mfso_4\oplus\mfso_4 &
\reals\oplus\reals\oplus\mfu_2 & \textrm{no} & \frac{4\pm\sqrt{6}}{5}& \frac{1}{X_1}\\

\xstrut\mfsp_{4l} & \mfsp_{2l}\oplus\mfsp_{2l} &
\mfsp_l\oplus\mfsp_l\oplus\mfsp_l\oplus\mfsp_l,\,l\geq 1 & \textrm{no} & \frac{4l+1\pm\sqrt{14l^2+7l+4}}{5(2l+1)}& \frac{l}{2l+1}.\frac{1}{X_1} \\

\xstrut\mfsp_{4l} & \mfsp_{2l}\oplus\mfsp_{2l} &
\mfsp_l\oplus\mfsp_l\oplus\mfu_{2l},\,l\geq 3 & \textrm{no} & \frac{2(4l+1)\pm\sqrt{4l^2-8l-1}}{5(2l+1)}& \frac{l}{2l+1}.\frac{1}{X_1}\\
\hline\ear$$

\footnotetext{See Theorem \ref{genIIall}. Compare with Tables \ref{eigIIexc}, \ref{eigIIclass} and \ref{bimII}.}
\end{minipage}\etab

\btab[h!]\caption{Einstein bisymmetric fibrations of Type II with $\mfg$ exceptional}\label{tabgenII}
\begin{minipage}{\textwidth}
$$\bar{|lll|c|c|c|} \hline\mfg & \mfk & \mfl & 4\textrm{-sym.} & X_1 & X_2\\ \hline
\xstrut\mfg_2 & \mfsu_2\oplus\mfsu_2 & \reals\oplus\reals & \textrm{no}  & 0.5526 & 3.6958\\
\xstrut & & & & 0.7432 & 4.7185\\\hline
 \xstrut \mfe_6 & \mfsu_6\oplus\mfsu_2 & \mfsu_5\oplus\reals\oplus\reals & \textrm{yes} & 1.5838 & 5.2195\\
\xstrut & & & & 0.3702 & 4.6215\\
\xstrut & & & & 0.5345 & 0.6682\\
\xstrut & & & & 1.0499 & 0.6338\\\hline
 \xstrut \mfe_7 & \mfso_{12}\oplus\mfsu_2 & \reals\oplus\mfso_{10}\oplus\reals & \textrm{yes} &  0.3086 & 7.4890 \\
\xstrut & & & & 0.4686 & 0.6737\\
\xstrut & & & & 0.9326 & 0.6496\\
\xstrut & & & & 1.4616 & 8.1878\\\hline
 \xstrut \mfe_7 & \mfso_{12}\oplus\mfsu_2& \mfso_{4}\oplus\mfso_{8}\oplus\reals & \textrm{no} & 0.3143 & 7.3931\\
\xstrut & & & & 1.4375 & 8.0839\\\hline
 \xstrut \mfe_7 & \mfso_{12}\oplus\mfsu_2& \mfso_{6}\oplus\mfso_{6}\oplus\reals & \textrm{yes} & 0.3163 & 7.3606\\
\xstrut & & & & 1.4292 & 8.0485\\  \hline\ear$$
\footnotetext{See Theorem \ref{genIIexc}. Compare with Table \ref{eigIIexc}.}
\end{minipage}\etab

\newpage
\section*{acknowledgements}
I would like to thank professor Dmitri Alekseevsky for his useful advice and enlightening discussions. This work was sponsored by Funda\c{c}\~{a}o para a Ci\^{e}ncia e a Tecnologia.

\end{document}